\documentclass[11pt,a4paper]{article}

\usepackage[margin=1in]{geometry}
\usepackage[T1]{fontenc}
\usepackage{lmodern}

\usepackage{amsmath,amssymb,amsthm}

\theoremstyle{plain}
\newtheorem{theorem}{Theorem}[section]
\newtheorem{lemma}[theorem]{Lemma}

\newtheorem{corollary}[theorem]{Corollary}

\theoremstyle{definition}
\newtheorem{definition}[theorem]{Definition}

\theoremstyle{remark}
\newtheorem{remark}[theorem]{Remark}

\newtheorem{assumption}{Assumption}

\usepackage{graphicx}

\usepackage{amsmath, amssymb, amsthm}
\usepackage{microtype}
\usepackage{graphicx}
\usepackage{hyperref}
\usepackage{xcolor}
\usepackage{amsfonts}

\usepackage{doi}
\usepackage{bbm}

\usepackage{subcaption}
\usepackage{caption}

\usepackage{algorithm}
\usepackage{algpseudocode}
\usepackage{float}

\usepackage{tikz}
\usetikzlibrary{positioning, calc}
\definecolor{myblue}{RGB}{31, 97, 158}
\definecolor{myred}{RGB}{192, 0, 0}
\definecolor{goldenscript}{RGB}{120, 50, 0}
\definecolor{mygreenfill}{RGB}{178, 224, 162}
\definecolor{mygreentext}{RGB}{30, 86, 27}
\definecolor{myorangefill}{RGB}{255, 230, 150}

\usepackage{comment}

\newcommand{\reals}{\mathbb{R}}

\newcommand{\Ftilde}{\tilde{F}}

\title{Direct Search meets Frank-Wolfe: Inexact Derivative-Free Bilevel Optimization over Polytopes}

\author{
Matteo Pernini\textsuperscript{*}
\quad Francesco Rinaldi\textsuperscript{\dag}
\quad Gabriele Sanguin\textsuperscript{\dag}
\quad Marco Viola\textsuperscript{\ddag}
}

\date{}

\begin{document}

\maketitle

\begingroup
\renewcommand{\thefootnote}{\fnsymbol{footnote}}

\footnotetext[1]{
ESSEC Business School, Paris, France.
E-mail: \texttt{matteo.pernini@essec.edu}.
}

\footnotetext[2]{
Department of Mathematics ``Tullio Levi-Civita'',
University of Padova, Padova, Italy.
E-mails: \texttt{rinaldi@math.unipd.it},
\texttt{gabriele.sanguin.2@phd.unipd.it}.
}

\footnotetext[3]{
School of Mathematical Sciences,
Dublin City University, Dublin, Ireland.
E-mail: \texttt{marco.viola@dcu.ie}.
}

\endgroup

\begin{abstract}
In this work we propose BLDS-SIMPLEX, a derivative-free direct search method for the solution of bilevel optimization with polytope constraints at the upper level. 
The proposed method reformulates the feasible set using barycentric coordinates and explores it using Frank-Wolfe-type directions, relying on a fixed-accuracy black-box oracle for the upper-level objective function. We analyze its convergence properties and  complexity in the case of smooth reduced objective; additionally, we analyze a variant  for the case of nonsmooth objectives and show its convergence to a specific kind of stationary point that we introduce as a natural generalization to the classic stationarity notions in this setting. Numerical tests on modified instances from  BOLIB support the viability of our new approach.

\end{abstract}


\bigskip

\section{Introduction}
\label{sec:introduction}

We consider bilevel optimization problems with polytope
constraints at the upper level:
\begin{equation}
\label{eq:blop}
\begin{aligned}
\min_{x \in \Omega_x} \quad & f(x,y^*(x)) \\
\text{s.t.} \quad &
y^*(x)\in\arg\min_{y\in\Omega_y}g(x,y),
\end{aligned}
\end{equation}
where $\Omega_x\subseteq\reals^{n_x}$ is a nonempty polytope
and $\Omega_y\subseteq\reals^{n_y}$ is the lower-level
feasible set. We focus on problems with a unique lower-level
response and access to approximate lower-level solutions
computed at a fixed accuracy. The upper-level algorithm
uses the resulting inexact objective evaluations without
requiring derivative information. The precise assumptions
on the problem and the oracle are stated in
Section~\ref{sec:preliminaries}.
  Our goal is to study the convergence properties of a derivative-free optimization (DFO) method in both smooth and nonsmooth settings.

\paragraph{Contributions.}
In this work, we make the following contributions:
\begin{enumerate}
    \item \textbf{Algorithm design.}
    We propose BLDS-SIMPLEX
    (Algorithm~\ref{alg:blds-simplex}), a derivative-free
    direct-search algorithm for bilevel optimization on polytopes.
    The method reformulates the upper-level problem on the unit
    simplex through an atomic representation of the feasible set
    and explores Frank-Wolfe-type directions using inexact
    upper-level function evaluations obtained through a
    lower-level oracle.

    \item \textbf{Convergence in the smooth setting.}
Under suitable smoothness assumptions and bounded
function-evaluation errors, we establish finite termination
when a positive minimum stepsize is enforced. More specifically, we derive
a bound on the Frank-Wolfe gap at the final iterate,
quantifying the effects of oracle accuracy and algorithmic
parameters
(Theorem~\ref{thm:smooth_fw_convergence_df}).
  \item \textbf{Analysis in the nonsmooth setting.}
We introduce the notion of a Nearby-Clarke Frank-Wolfe
stationary point (NCFWSP) and study convergence properties
of a suitable variant of the method using this criterion
(Theorem~\ref{thm:gfw_convergence}).
The analysis builds on the inexact direct-search framework
of \cite{diouane2024inexact}, accounting for oracle
inexactness and polytope constraints at the upper level.

    \item \textbf{Numerical experiments.}
    We evaluate the proposed method 
    on modified bilevel benchmark problems from the BOLIB
    library \cite{bolib3,zhou2020bolib} and compare its
    empirical performance with that of existing derivative-free
    approaches.
\end{enumerate}

\paragraph{Paper outline.}
The rest of the paper is organized as follows. Section~\ref{sec:related_work} reviews the related literature. Section~\ref{sec:preliminaries} introduces the problem formulation, the algorithm, and some necessary assumptions and definitions. Section~\ref{sec:theoretical_analysis} presents the convergence analysis in both the smooth and nonsmooth settings. Section~\ref{sec:numerical_experiments} reports the numerical results. Section~\ref{sec:conclusion} concludes the paper.

\section{Related Work}
\label{sec:related_work}

\paragraph{Bilevel optimization.}
Bilevel optimization has a long history in operations research and mathematical programming \cite{colson2007overview,dempe2002foundations}. 
Interest in this topic has grown significantly in recent years due to applications in machine learning, such as hyperparameter optimization and meta-learning (see, e.g., \cite{franceschi2017forward,franceschi2018bilevel} and references therein).
For differentiable convex lower-level problems satisfying suitable
constraint qualifications, classical approaches replace the
lower-level optimality requirement with its KKT conditions,
yielding an equivalent mathematical program with equilibrium
 constraints (MPEC) \cite{colson2007overview}. 
More recently, gradient-based methods that exploit explicit or
implicit differentiation to compute hypergradients have gained
popularity in machine learning (see, e.g.,
\cite{zhang2024introduction} and references therein for a survey).
Many implicit-differentiation approaches require suitable
differentiability assumptions and access to derivative information. These derivative requirements
can limit applicability in black-box settings, where such
information may be unavailable or costly to obtain. Moreover,
hypergradient calculation can be computationally expensive and
sensitive to the conditioning of the lower-level problem.
For these reasons, our work studies a derivative-free method that only requires inexact function evaluations of the upper-level objective, making it applicable to a broader class of bilevel problems with black-box structure.

\paragraph{Derivative-free methods.}
Derivative-free optimization (DFO) methods are designed to minimize functions without requiring or using derivative information. They can be broadly categorized into direct-search, model-based and randomized methods (see, e.g., \cite{audet2014survey,conn2009introduction,dzahini2025direct,larson2019derivative} and references therein for further details on these three classes of methods).

In the context of bilevel optimization, several derivative-free approaches have been developed. To cite some, Mersha and Dempe \cite{mersha2011direct} and Zhang and Lin \cite{zhang2014bilevel} applied direct-search methods to bilevel problems with various assumptions on the lower-level oracle, respectively assuming no approximation error and adaptive error.
Conn and Vicente \cite{conn2012bilevel} developed a trust-region
framework accommodating inexact lower-level solutions and
reuse of sample points. Ehrhardt and Roberts
\cite{ehrhardt2021inexact} analyzed a dynamically accurate
derivative-free trust-region method and established
worst-case complexity bounds.
More recently, zeroth-order methods based on smoothing strategies \cite{chen2023bilevel,maheshwari2023convergent} have been studied, offering tailored complexity estimates.
Finally, min-max DFO problems \cite{anagnostidis2021direct,menickelly2020derivative} and direct-search methods under adaptive deterministic noise \cite{lucidi2002derivative} and stochastic noise \cite{audet2019stomads,rinaldi2024stochastic} have also been tackled in the literature.

Most closely related to our work is the analysis of
Diouane et al.~\cite{diouane2024inexact}, who developed
directional and mesh-based direct-search methods for bilevel
problems with a fixed-accuracy lower-level oracle.
They analyzed convergence in both smooth and nonsmooth
settings and provided complexity estimates. Their analysis, however, concerns unconstrained
problems. Our work builds on this framework by considering
polytope constraints at the upper level while retaining
the fixed-accuracy oracle model.

To the best of our knowledge, previous analyses of direct-search methods for bilevel optimization have not simultaneously considered bounded, non-adaptive lower-level inexactness and polytope constraints at the upper level. Our approach differs from existing methods by explicitly exploiting the polytope geometry through Frank-Wolfe-type directions. Consequently, its convergence guarantees are expressed in terms of Frank-Wolfe gap-based stationarity measures. 
In this respect, our approach shares some similarities with the simplicial direct-search framework introduced by Cristofari and Rinaldi~\cite{cristofari2021derivative}. However, their scheme exploits barycentric coordinates and sparse feasible directions to guarantee convergence for single-level black-box problems, while we extend these concepts to an inexact bilevel setting where approximate lower-level solutions lead to inexact upper-level evaluations.

In a recent work, Cesaroni et al.~\cite{cesaroni2026derivativefreebileveloptimizationinexact} proposed a derivative-free framework for bilevel optimization that accounts for lower-level inexactness. Their framework differs from ours in two main respects. First, our algorithm is a feasible method designed to directly handle upper-level feasible regions given by general polytopes, whereas Cesaroni et al.\ focus primarily on bound-constrained problems and treat general nonlinear upper-level constraints via an exact-penalty reformulation. Second, their analysis is built around a progressive refinement scheme for lower-level accuracy, establishing convergence to Clarke-Jahn stationary points as the lower-level error tends to zero and approximate Goldstein stationarity when refinement is limited, whereas our approach directly analyzes a fixed  lower-level inexactness.

\section{Preliminaries}
\label{sec:preliminaries}

\subsection{Problem Formulation}
\label{subsec:problem_formulation}
We now address the general formulation of the bilevel optimization problem~\eqref{eq:blop}, with a nonempty compact upper-level feasible set $\Omega_x$.
We will later specialize to the case where $\Omega_x$ is a polytope described as the convex hull of a finite set of atoms.
First, notice that there might be different solutions to the lower-level problem, which can lead to different upper-level objectives. To ensure that the upper-level objective is well-defined and to simplify the analysis, it is common practice in the literature to assume that the lower-level problem has a unique solution for each $x$, which is known as the lower-level singleton (LLS) assumption (see, e.g., \cite{diouane2024inexact} and references therein for further details).

\begin{assumption}[LLS]
\label{ass:lls}
For any $x \in \Omega_x$, we have that $\arg\min_{y \in \Omega_y} g(x, y) = \{y^*(x)\} $.
\end{assumption}

\noindent Under this assumption, problem~\eqref{eq:blop} reduces to
\begin{equation}
\label{eq:blop_reduced}
\min_{x \in \Omega_x} F(x) := f(x, y^*(x)).
\end{equation}

In practical implementations, it is often not possible to compute the lower-level solution $y^*(x)$ exactly. Instead, we assume to have access to an approximate solution $\tilde{y}(x)$ leading to the inexact objective $\Ftilde(x) := f(x, \tilde{y}(x))$. By calling this \emph{lower-level oracle}, the algorithm always operates on $\Ftilde$, bringing the need for the following two assumptions to analyze its convergence in terms of the true objective $F$:

\begin{assumption}[Bounded Lower-Level Error]
\label{ass:bounded_lower_level_error}
For any $x \in \Omega_x$, the approximate lower-level solution $\tilde{y}(x)\in \Omega_y$ satisfies $\|\tilde{y}(x) - y^*(x)\| \le \varepsilon$ for a given $\varepsilon > 0$.
\end{assumption}

\begin{assumption}[Lipschitz Continuity of $f$ w.r.t.\ $y$]
\label{ass:lipschitz_f_y}
The function $f(x, y)$ is $L_f$-Lipschitz continuous with respect to $y$, uniformly in $x$:
\[
|f(x, y_1) - f(x, y_2)| \leq L_f \|y_1 - y_2\|, \quad \forall\, y_1, y_2 \in \Omega_y,\; \forall\, x \in \Omega_x.
\]
\end{assumption}

\noindent Assumption~\ref{ass:bounded_lower_level_error} says that the lower-level oracle always gives an approximate solution within a given error bound $\varepsilon$ from the exact solution. We remark that this can be satisfied by applying standard solvers to the lower-level problem, e.g. gradient descent or DFO methods with suitable error bounds (see discussion in \cite{diouane2024inexact}).
With Assumption~\ref{ass:lipschitz_f_y}, we are giving a bound on the change in the upper-level objective function due to changes in the lower-level solution. By combining the two, we get an error bound between the inexact and exact objectives:

\begin{remark}[Bounded objective error]
\label{rem:bounded_objective_error}
Assumptions~\ref{ass:bounded_lower_level_error} and~\ref{ass:lipschitz_f_y} imply: 
\[ |\tilde{F}(x) - F(x)| \le L_f \varepsilon \text{ for all } x \in \Omega_x \].
\end{remark}

To study convergence both in the smooth and nonsmooth case, we also need the following continuity assumption on the upper-level objective function $F$.

\begin{assumption}[Lipschitz Continuity of $F$]
\label{ass:lipschitz_f_x}
$F(x)$ is $L_F$-Lipschitz continuous in $\Omega_x$, i.e.
    \[|F(x_1) - F(x_2)| \leq L_F \|x_1 - x_2\|, \quad \forall x_1, x_2 \in \Omega_x.\]
\end{assumption}

\noindent Note that if $f$ is Lipschitz continuous with respect to $x$, and $y^*(x)$ is Lipschitz continuous with respect to $x$, then Assumption~\ref{ass:lipschitz_f_x} is satisfied.

\begin{remark}[Lower boundedness of $F$]
\label{rem:f_min}
Under Assumption~\ref{ass:lipschitz_f_x}, the objective function $F(x)$ is bounded below on $\Omega_x$ by some constant $f_{\min} > -\infty$, due to the compactness of the domain (we assume $\Omega_x$ closed and bounded) and the continuity of $F$.
By the error bound of Remark~\ref{rem:bounded_objective_error}, we also have that
\[
\tilde{F}(x) \ge \tilde{f}_{\min} := f_{\min} - L_f \varepsilon > -\infty \text{ for all } x \in \Omega_x.
\]
\end{remark}

As anticipated above, from now on we are considering the upper-level feasible set $\Omega_x \subseteq \reals^{n_x}$ to be a polytope described by a finite set of \emph{atoms} $\mathcal{A} = \{a_1, \ldots, a_m\} \subset \reals^{n_x}$, so that $\Omega_x = \mathrm{conv}(\mathcal{A})$.
Any point $x \in \Omega_x$ can therefore be described as a convex combination $x = Aw$, where $A = [a_1 \; a_2 \; \cdots \; a_m] \in \reals^{n_x \times m}$ is the atom matrix and $w \in \Delta^{m-1} := \{w \in \reals^m : w \ge 0,\; \mathbf{1}^\top w = 1\}$ is the vector of barycentric coordinates. Thus, problem~\eqref{eq:blop} can be equivalently reformulated considering the simplicial representation of the feasible set:
\begin{equation}
\label{eq:blop_simplex}
\begin{array}{rl}
    \underset{w \in \Delta^{m-1}}{\min} & f(Aw, y^*(Aw)) \\
    \text{s.t.} & y^*(Aw) \in \underset{y \in \Omega_y}{\arg \min} \; g(Aw, y).
\end{array}
\end{equation}

\noindent Equivalently, we can write the reduced problem~\eqref{eq:blop_reduced} as

\begin{equation}
\label{eq:blop_simplex_reduced}
\min_{w \in \Delta^{m-1}} F_A(w) := f(Aw, y^*(Aw)).
\end{equation}

The assumptions above transfer to the simplex formulation
through the map $x=Aw$. In particular, the inexact objective is
\[
\Ftilde_A(w):=\Ftilde(Aw)=f(Aw,\tilde y(Aw)).
\]
By Remark~\ref{rem:bounded_objective_error}, the objective-error
bound is preserved:
\[
|\Ftilde_A(w)-F_A(w)|\le L_f\varepsilon,
\qquad \forall w\in\Delta^{m-1}.
\]
Moreover, Assumption~\ref{ass:lipschitz_f_x} gives
\[
|F_A(w_1)-F_A(w_2)|
\le L_F\|A\|\,\|w_1-w_2\|,
\qquad \forall w_1,w_2\in\Delta^{m-1}.
\]
Thus, $L_{F_A}:=L_F\|A\|$ is a valid Lipschitz constant for
$F_A$. Here vector norms are Euclidean and $\|A\|$ denotes
the induced operator norm.

\subsection{Algorithm}
\label{subsec:algorithm}

The proposed algorithm, BLDS-SIMPLEX, is a \emph{BiLevel Direct-Search method on the Unit Simplex} that iteratively explores Frank-Wolfe-type directions $d = \pm(e_i - w^k)$, where $e_i$ is the $i$-th standard basis vector. At each iteration, the algorithm cycles through a permutation of the coordinate indices and, for each index, tests both positive and negative directions.

 The pseudocode of the proposed method is reported in Algorithm~\ref{alg:blds-simplex}. We recall the definition of the maximum stepsize that guarantees feasibility:
\begin{definition}[Maximum Feasible Step]
    \label{def:max_feasible_step}
    Let $w^k\in \Delta^{m-1}$. We define the maximum feasible step for $d \in \reals^m$, as
    
    \begin{equation}
    \bar{\alpha}(w^k,d) := \max\{\alpha \geq 0: w^k + \alpha d \in \Delta^{m-1}\}.
    \end{equation}
Note that for the specific case of FW directions in the simplex, the computation of the maximal feasible stepsize has the closed form:
\[\bar{\alpha}(w^k, d) = 
\begin{cases} 
    1 & \text{if } d = e_i - w^k,\\
    \frac{w^k_i}{1 - w^k_i} & \text{if } d = -(e_i - w^k).
\end{cases}\]
\end{definition}
In BLDS-SIMPLEX, a trial point $w^k+\alpha d$ is accepted if it produces a sufficient decrease in the inexact objective:
\[
\Ftilde_A(w^k+\alpha d)
\le \Ftilde_A(w^k)-\gamma(\alpha^k)^2,
\]
where $\gamma>0$ is a forcing parameter and
$\alpha=\min\{\alpha^k,\bar\alpha(w^k,d)\}$ is the effective stepsize, truncated at the maximum feasible step $\bar\alpha(w^k,d)$ to keep the trial point on the simplex.
Upon a successful step, the nominal stepsize is enlarged, subject to the prescribed maximum; if all polling directions fail, it is contracted, subject to the minimum stepsize.
When a positive minimum stepsize $\alpha_{\min}>0$ is enforced, the algorithm terminates after an unsuccessful iteration with $\alpha^k=\alpha_{\min}$, at which point the stepsize can no longer be reduced.

\begin{algorithm}[htbp]
\caption{BLDS-SIMPLEX} 
\label{alg:blds-simplex}
\begin{algorithmic}[1]
\Require Choose $w^0 \in \Delta^{m-1}$, initial stepsize $\alpha^0 \in (0,1]$, minimum stepsize $\alpha_{\min} \in [0,1]$, maximum stepsize $\alpha_{\max} \in (0,1]$, parameters $\theta \in (0,1)$, $\delta \in (0,1)$, and forcing constant $\gamma>0$
\For{$k = 0,1,2,\ldots$}
    \State Set $\texttt{success} = \texttt{false}$ and $\alpha^k = \max\{\alpha^k, \alpha_{\min}\}$
    \State Choose a permutation of $\{1, \ldots, m\}$ and denote its elements as $\{i_j\}_{j=1}^{m}$. Set $j=1$
    \While {$j \le m$}
        \For{$\mathrm{sgn} \in \{+, -\}$}
            \State Set $d = (\mathrm{sgn})(e_{i_j} - w^k)$
            \State Compute $\bar{\alpha}(w^k, d)$, the largest $\alpha$ such that $w^k + \alpha d \in \Delta^{m-1}$
            \State Set $\hat{\alpha} = \min\{\alpha^k, \bar{\alpha}(w^k, d)\}$
            \If{$\hat{\alpha} > 0$}
                \State Call lower-level oracle to compute $\tilde{y}(A(w^k + \hat{\alpha} d))$ and evaluate $\Ftilde_A(w^k + \hat{\alpha} d)$
                \If{$\Ftilde_A(w^k + \hat{\alpha} d) \le \Ftilde_A(w^k) - \gamma (\alpha^k)^2$}
                    \State Set $w^{k+1} = w^k + \hat{\alpha} d$ and $\texttt{success} = \texttt{true}$
                    \State \textbf{break} (exit the sign loop)
                \EndIf
            \EndIf
        \EndFor
        \If{$\texttt{success}$}
            \State \textbf{break} (exit the while loop)
        \EndIf
        \State Set $j = j + 1$
    \EndWhile
    \If{$\texttt{success}$}
        \State Set $\alpha^{k+1} = \min\{\alpha_{\max}, \alpha^k/\delta\}$ \Comment{Expand stepsize}
    \Else
        \State Set $w^{k+1} = w^k$ and $\alpha^{k+1} = \max\{\alpha_{\min}, \theta \alpha^k\}$ \Comment{Shrink stepsize}
        \If{some approximate stationarity condition satisfied}
            \State \textbf{Terminate} \Comment{Stopping criterion}
        \EndIf
    \EndIf
\EndFor
\end{algorithmic}
\end{algorithm}

\section{Theoretical Analysis}
\label{sec:theoretical_analysis}

In this section, we present the convergence analysis of Algorithm~\ref{alg:blds-simplex}. We first consider the smooth setting (Section~\ref{subsec:smooth_analysis}), where the reduced objective $F$ has a Lipschitz continuous gradient, and then the nonsmooth setting (Section~\ref{subsec:nonsmooth_analysis}), where only Lipschitz continuity of $F$ is assumed.

\subsection{Convergence in the Smooth Setting}
\label{subsec:smooth_analysis}

In this subsection, we impose the following smoothness
assumption on the reduced objective.

\begin{assumption}[Smoothness of $F$]
\label{ass:smooth_F}
The reduced objective $F(x)=f(x,y^*(x))$ admits a continuously
differentiable extension to an open neighborhood of $\Omega_x$,
denoted again by $F$. Its gradient is $L$-Lipschitz continuous
on $\Omega_x$, i.e.,
\[
\|\nabla F(x_1)-\nabla F(x_2)\|
\le L\|x_1-x_2\|,
\qquad \forall x_1,x_2\in\Omega_x.
\]
\end{assumption}

For the simplex reformulation $F_A(w)=F(Aw)$, the chain rule
gives
\[
\nabla F_A(w)=A^\top\nabla F(Aw).
\]
Consequently, $F_A$ has a Lipschitz continuous gradient on
$\Delta^{m-1}$ with a valid constant
\[
L_A:=L\|A\|^2.
\]

This smoothness assumption allows us to apply the descent
lemma to control directional derivatives along the forward
Frank--Wolfe directions at unsuccessful iterations.
These bounds yield a stationarity guarantee expressed
in terms of the Frank--Wolfe gap, which we recall next.

\begin{definition}[Frank-Wolfe Gap]
\label{def:fw_gap}
For a differentiable function $F$ on a convex set $\Omega$, the Frank-Wolfe gap at $x \in \Omega$ is
\begin{equation}
    \label{eq:fw_gap}
g_{\mathrm{FW}}(x; F, \Omega) := \max_{u \in \Omega} \langle -\nabla F(x),\, u - x \rangle.
\end{equation}
A point $x^*$ is $\epsilon$-Frank-Wolfe stationary if $g_{\mathrm{FW}}(x^*; F, \Omega) \le \epsilon$.
\end{definition}

\noindent The Frank-Wolfe gap is a common criterion for measuring convergence rate of projection-free algorithms in the smooth setting: it vanishes at first-order stationary points and can be efficiently computed for polytopes, as it reduces to a linear optimization problem over the feasible set.

Let us start now with the convergence analysis of Algorithm~\ref{alg:blds-simplex}. To recall everything seen until now, we are considering the reduced bilevel problem~\eqref{eq:blop_simplex_reduced} under Assumptions~\ref{ass:lls},~\ref{ass:bounded_lower_level_error},~\ref{ass:lipschitz_f_y} and \ref{ass:smooth_F}, which will be referred as the \emph{smooth setting assumptions}.
Note that Assumption~\ref{ass:smooth_F} directly implies Assumption~\ref{ass:lipschitz_f_x} on the compact domain $\Omega_x$. 
Our goal is to show convergence to a point with a small Frank-Wolfe gap~(\ref{def:fw_gap}) and give the complexity bounds when $\alpha_{\min} \geq 0$.
We begin with a lemma that bounds the steepest feasible \emph{FW-direction} at an unsuccessful iteration.

\begin{lemma}[Gradient bound at unsuccessful iterations]
\label{lem:fw_grad_bound}
Let the smooth setting assumptions hold, and let
$\mathcal{D}^k_\Delta
:=\{e_i-w^k:i=1,\ldots,m\}
$
be the set of elementary forward FW-directions at iteration
$k$. If iteration $k$ is unsuccessful, then, for every
$d\in\mathcal{D}^k_\Delta$,
\begin{equation}
\label{eq:fw_grad_bound}
-\nabla F_A(w^k)^\top d
\le (L_A+\gamma)\alpha^k
+\frac{2L_f\varepsilon}{\alpha^k}.
\end{equation}
\end{lemma}

\begin{proof}
At iteration $k$, consider any FW-direction
$d=e_i-w^k\in\mathcal{D}^k_\Delta$.
If $d=0$, the claimed inequality holds immediately.
Otherwise, the maximum feasible stepsize along $d$ is one.
Since $0<\alpha^k\le1$, the algorithm therefore evaluates
the feasible trial point $w^k+\alpha^k d$ without truncating
the nominal stepsize.

Because the iteration is unsuccessful, we have
\begin{equation}
\label{eq:unsuccessful_iter}
\Ftilde_A(w^k+\alpha^k d)
>
\Ftilde_A(w^k)-\gamma(\alpha^k)^2.
\end{equation}

From Remark~\ref{rem:bounded_objective_error}, we have
\begin{align*}
\Ftilde_A(w^k)
&\ge F_A(w^k)-L_f\varepsilon,\\
\Ftilde_A(w^k+\alpha^k d)
&\le F_A(w^k+\alpha^k d)+L_f\varepsilon.
\end{align*}
Combining these inequalities with
\eqref{eq:unsuccessful_iter}, we obtain
\begin{equation}
\label{eq:true_diff_lower}
F_A(w^k+\alpha^k d)-F_A(w^k)
>
-\gamma(\alpha^k)^2-2L_f\varepsilon.
\end{equation}

By Assumption~\ref{ass:smooth_F} and the linear
reformulation, $F_A$ has an $L_A$-Lipschitz continuous
gradient on $\Delta^{m-1}$. Thus, the descent lemma gives
\begin{equation}
\label{eq:descent_lemma}
\begin{aligned}
F_A(w^k+\alpha^k d)-F_A(w^k)
&\le \alpha^k\nabla F_A(w^k)^\top d\\
&\quad+\frac{L_A}{2}(\alpha^k)^2\|d\|^2.
\end{aligned}
\end{equation}

Combining \eqref{eq:true_diff_lower} and
\eqref{eq:descent_lemma}, and dividing by $\alpha^k>0$,
we find
\[
-\nabla F_A(w^k)^\top d
<
\left(\gamma+\frac{L_A}{2}\|d\|^2\right)\alpha^k
+\frac{2L_f\varepsilon}{\alpha^k}.
\]

Finally, since $w^k\in\Delta^{m-1}$, we have
$\|d\|^2\leq 2.$
Consequently,
\[
-\nabla F_A(w^k)^\top d
\le (L_A+\gamma)\alpha^k
+\frac{2L_f\varepsilon}{\alpha^k},
\]
which proves the claimed bound.
\end{proof}

We now state the main convergence result for the smooth case.

\begin{theorem}[Smooth convergence and FW-gap bound]
\label{thm:smooth_fw_convergence_df}
Let the smooth setting assumptions hold. Then:
\begin{enumerate}
    \item If $\alpha_{\min} > 0$, Algorithm~\ref{alg:blds-simplex}
    terminates at an iteration $\bar{k}$ satisfying
    \[
    \bar{k} \le 1 + \frac{\Ftilde_A(w^0) - \tilde{f}_{\min}}{\gamma \alpha_{\min}^2} \left(2 + \frac{\ln \delta}{\ln \theta}\right) + \frac{\ln \alpha_{\min} - \ln \alpha^0}{\ln \theta},
    \]
    and the final iterate $w^{\bar{k}}$ satisfies:
    \begin{equation}
    \label{eq:fw_gap_smooth}
    g_{\mathrm{FW}}(w^{\bar{k}}; F_A, \Delta^{m-1}) \le \Big(L_A + \gamma\Big) \alpha_{\min} + \frac{2L_f \varepsilon}{\alpha_{\min}}.
    \end{equation}

    \item If $\alpha_{\min} = 0$, then $\alpha^k \to 0$.
    Moreover, for any $\bar{\alpha} > 0$ with
    $\alpha^0 \ge \bar{\alpha}$, there exists an unsuccessful
    iteration $\bar{k}$ such that:
    \begin{equation}
    \label{eq:fw_gap_alpha0}
    g_{\mathrm{FW}}(w^{\bar{k}}; F_A, \Delta^{m-1}) < \frac{1}{\theta}\Big[\Big(L_A + \gamma\Big) \bar{\alpha} + \frac{2L_f \varepsilon}{\bar{\alpha}}\Big],
    \end{equation}
    and for all $k \ge \bar{k}$,
    \[
    F_A(w^k) \le F_A(w^{\bar{k}}) + 2L_f \varepsilon.
    \]
\end{enumerate}
\end{theorem}

\begin{proof}
    Consider the sequence $\{\tilde{F}_A(w^k)\}$. It is
    nonincreasing by definition and bounded below
    (see Remark~\ref{rem:f_min}); therefore, whenever the
    algorithm generates an infinite sequence, these values
    converge to some finite limit $\ell \geq \tilde{f}_{\min}$.

    Suppose now $\alpha_{\min} > 0$. It follows that the
    number $k_s$ of successful iterations is finite:
    \begin{equation}
        \label{eq:ks_bound}
    k_s \leq \frac{\tilde{F}_A(w^0) - \tilde{f}_{\min}}{\gamma \alpha_{\min}^2},
    \end{equation}
    where we used the fact that, for each successful step,
    the decrease in the objective is at least
    $\gamma \alpha_{\min}^2$.
    After the last successful iteration, consecutive
    unsuccessful iterations reduce the stepsize to
    $\alpha_{\min}$, and the next unsuccessful iteration
    terminates the algorithm. Thus termination occurs
    after finitely many iterations.

    Let $k_{ns}$ be the number of unsuccessful iterations,
    so that $k_s + k_{ns} = \bar{k}+1$ is the total number
    of iterations. Among the nonterminal unsuccessful
    iterations, let $k_f$ count those for which
    $\theta\alpha^k \geq \alpha_{\min}$, and let $k_c$
    count those for which $\theta\alpha^k < \alpha_{\min}$.

    On a successful iteration, the stepsize is increased
    by at most a factor $1/\delta$; on each of the $k_f$
    iterations, it is multiplied by $\theta$; and the
    remaining updates cannot increase it. Since the final
    stepsize is $\alpha_{\min}$, we have
    \begin{equation}
        \label{eq:alpha_bound}
        \alpha_{\min} \leq \frac{\theta^{k_f}}{\delta^{k_s}} \alpha^0.
    \end{equation}
    Moreover, each of the $k_c$ clipped contractions sets
    the next stepsize to $\alpha_{\min}$. Before another
    such contraction can occur, there must be a successful
    iteration; otherwise, the next unsuccessful iteration
    would terminate the algorithm. Hence $k_c \leq k_s+1$.

    Including the terminal unsuccessful iteration, we have
    $k_{ns}=k_f+k_c+1$. Using \eqref{eq:alpha_bound}, we
    can therefore bound $k_{ns}$ as follows:
    \begin{align}
        k_{ns} &\leq 2 + k_s + \frac{1}{\ln(\theta)} \left(\ln(\alpha_{\min}) - \ln(\alpha^0) + k_s \ln(\delta)\right) \nonumber\\
        &\leq 2 + \frac{\tilde{F}_A(w^0)-\tilde{f}_{\min}}{\gamma\alpha_{\min}^2}
        \left(1+\frac{\ln(\delta)}{\ln(\theta)}\right)
        + \frac{\ln(\alpha_{\min})-\ln(\alpha^0)}{\ln(\theta)}.
        \label{eq:kns_bound}
    \end{align}

    We then combine the bounds on successful and
    unsuccessful iterations \eqref{eq:ks_bound} and
    \eqref{eq:kns_bound} to obtain the desired bound
    on $\bar{k}$:
    \begin{align}
        \bar{k} &= k_s + k_{ns}-1 \nonumber\\
        &\leq 1 + \frac{1}{\gamma \alpha_{\min}^2} (\tilde{F}_A(w^0) - \tilde{f}_{\min}) \left(2 + \frac{\ln(\delta)}{\ln(\theta)}\right) + \frac{\ln(\alpha_{\min}) - \ln(\alpha^0)}{\ln(\theta)}.
        \label{eq:k_bound}
    \end{align}

    To prove the bound on the FW-gap, we use
    Lemma~\ref{lem:fw_grad_bound}, which holds for every
    unsuccessful iteration. From the algorithm stopping
    criterion, the last iteration $\bar{k}$ is unsuccessful
    with $\alpha^{\bar{k}} = \alpha_{\min}$. Therefore,
    for every forward FW-direction
    $d \in \mathcal{D}_{\Delta}^{\bar{k}}$,
    Equation~\eqref{eq:fw_grad_bound} gives:
    \[
    -\nabla F_A(w^{\bar{k}})^\top d \;\le\; \Big(L_A+\gamma\Big)\,\alpha_{\min} \; + \; \frac{2L_f\,\varepsilon}{\alpha_{\min}}.
    \]

    Suppose now there exists at least one point
    $\hat{u} \in \Delta^{m-1}$ such that
    $$\big\langle -\nabla F_A(w^{\bar{k}}),\, \hat{u} - w^{\bar{k}}\big\rangle > 0.$$
    If no such point exists, we immediately obtain
    $$g_{\mathrm{FW}}(w^{\bar{k}};F_A,\Delta^{m-1}) = 0,$$
    and the claimed bound follows. Otherwise,
    \[
    \begin{aligned}
    g_{\mathrm{FW}}(w^{\bar{k}};F_A,\Delta^{m-1})
    &:= \max_{u\in\Delta^{m-1}} \;\big\langle -\nabla F_A(w^{\bar{k}}),\, u - w^{\bar{k}}\big\rangle \\
    &= \max_{u\in\Delta^{m-1}} \;\Big\langle -\nabla F_A(w^{\bar{k}}),\, \sum_{i=1}^m u_i \,(e_i - w^{\bar{k}})\Big\rangle
    \\
    &= \max_{u\in\Delta^{m-1}} \;\sum_{i=1}^m u_i \,\big\langle -\nabla F_A(w^{\bar{k}}),\, e_i - w^{\bar{k}}\big\rangle \\
    &\le \max_{u\in\Delta^{m-1}} \;\sum_{i=1}^m u_i \, \max_{j\in\{1,\dots,m\}} \big\langle -\nabla F_A(w^{\bar{k}}),\, e_j - w^{\bar{k}}\big\rangle \\
    &\le \max_{j\in\{1,\dots,m\}} \;\big\langle -\nabla F_A(w^{\bar{k}}),\, e_j - w^{\bar{k}}\big\rangle
    \\[1mm]
    &\le \max_{d \in \mathcal{D}_{\Delta}^{\bar{k}}} \; \big(-\nabla F_A(w^{\bar{k}})^\top d\big).
    \end{aligned}
    \]

    Combining the previous inequalities, we obtain
    the FW-gap bound:
    \[
    g_{\mathrm{FW}}(w^{\bar{k}};F_A,\Delta^{m-1})
    \;\le\; \Big(L_A+\gamma\Big)\,\alpha_{\min} \; + \; \frac{2L_f\,\varepsilon}{\alpha_{\min}}.
    \]
    This concludes the proof for the case
    $\alpha_{\min} > 0$.

    Let us now consider the case $\alpha_{\min} = 0$.
    In this case, the algorithm does not terminate in
    a finite number of iterations. Summing the sufficient
    decreases over the successful iterations gives
    \[
    \gamma\sum_{k\in\mathcal S}(\alpha^k)^2
    \leq \tilde F_A(w^0)-\tilde f_{\min},
    \]
    where $\mathcal S$ is the set of successful iteration
    indices. Hence, for every $\tau>0$, there are only
    finitely many successful iterations with
    $\alpha^k\geq\delta\tau$. Beyond the last such iteration,
    every successful iteration produces
    $\alpha^{k+1}<\tau$, while unsuccessful iterations
    contract the stepsize. Consequently, the stepsize
    eventually falls below $\tau$ and remains below it.
    Since $\tau>0$ is arbitrary, $\alpha^k\to0$.

    Assuming $\alpha^0 \geq \bar{\alpha}$, let $j\geq1$
    be the first index such that
    $\alpha^j<\theta\bar{\alpha}$.
    Since successful iterations cannot decrease the
    stepsize, iteration $j-1$ is unsuccessful and
    $\alpha^j=\theta\alpha^{j-1}$. Setting
    $\bar{k}=j-1$, we therefore have
    $\alpha^{\bar{k}}\in[\theta\bar{\alpha},\bar{\alpha})$.

    Define now $\phi(\alpha)=a\alpha+b/\alpha$, with
    $a=L_A+\gamma$ and $b=2L_f\varepsilon$;
    $\phi$ is continuous and convex on $\mathbb{R}_{>0}$.
    Note that, for $\theta\in(0,1)$,
    \[
    \phi(\theta \bar{\alpha}) < \frac{1}{\theta}\phi(\bar{\alpha}) \qquad \text{and} \qquad \phi(\bar{\alpha}) < \frac{1}{\theta}\phi(\bar{\alpha}).
    \]
    Therefore,
    $\phi(\alpha^{\bar{k}})
    \leq \max\{\phi(\bar{\alpha}),\phi(\theta\bar{\alpha})\}
    < \frac{1}{\theta}\phi(\bar{\alpha})$.
    Applying Lemma~\ref{lem:fw_grad_bound} at the
    unsuccessful iteration $\bar{k}$ and using the
    simplex representation of the FW-gap above gives
    Equation~\eqref{eq:fw_gap_alpha0}.

    Finally, for any $\bar{k}\in\mathbb{N}_0$ and
    $k\geq\bar{k}$, we have:
    \[
    F_A(w^k) \leq \tilde{F}_A(w^k) + L_f \varepsilon \leq \tilde{F}_A(w^{\bar{k}}) + L_f \varepsilon \leq F_A(w^{\bar{k}}) + 2 L_f \varepsilon.
    \]
\end{proof}

\begin{corollary}[Optimal stepsize choice]
\label{cor:optimal_alpha}
The FW-gap bound in~\eqref{eq:fw_gap_smooth} is minimized
by choosing
\[
\alpha_{\min}^* = \sqrt{\frac{2L_f \varepsilon}{L_A + \gamma}},
\]
provided $L_f\varepsilon>0$ and
$\alpha_{\min}^*\leq\alpha^0$.
This choice yields:
\[
g_{\mathrm{FW}}(w^{\bar{k}}; F_A, \Delta^{m-1}) \le 2\sqrt{2L_f \varepsilon \Big(L_A + \gamma\Big)} = O(\sqrt{\varepsilon}),
\]
as $\varepsilon\downarrow0$ with the other constants fixed.
\end{corollary}

\subsection{Convergence in the Nonsmooth Setting}
\label{subsec:nonsmooth_analysis}
In this subsection, we relax Assumption~\ref{ass:smooth_F}
and require only Lipschitz continuity of the reduced
objective $F$, together with
Assumptions~\ref{ass:lls}--\ref{ass:lipschitz_f_x}.
This case leads to a nonsmooth polytope-constrained upper-level problem in which we do not have access to an exact lower-level solution. To adapt to this scenario, we need to make some changes to the original BLDS-SIMPLEX algorithm and to introduce an appropriate measure of stationarity.

Without differentiability, the Frank-Wolfe gap is no longer well-defined, because of its dependence on the gradient. Therefore, we introduce the concept of $(\delta,\epsilon)$-Nearby-Clarke Frank-Wolfe Stationary Point (NCFWSP), aimed at extending the Frank-Wolfe gap to nonsmooth objectives through the Clarke subdifferential. We first recall some definitions.

\begin{definition}[Clarke Subdifferential \cite{clarke1990optimization}]
\label{def:clarke}
Let $F:  \reals^n \to \reals$ be locally Lipschitz at $x$. The Clarke subdifferential of $F$ at $x$ is
\[
\partial F(x) := \{g \in \reals^n : DF(x; v) \ge \langle g, v \rangle \;\; \forall\, v \in \reals^n\},
\]
where $D F(x; v) := \limsup_{y \to x,\, t \downarrow 0} \frac{F(y + tv) - F(y)}{t}$ is the generalized directional derivative. Each element $g$ of $\partial F(x)$ is called a generalized gradient of $F$ at $x$.
\end{definition}

\begin{definition}[Clarke Frank-Wolfe gap \cite{liu2024zeroth}]
\label{def:CFWgap}

Let $\Omega\subseteq\reals^n$ be nonempty, compact, and
convex, and let $F$ be locally Lipschitz on an open
neighborhood of $\Omega$. For $x\in\Omega$, the Clarke
Frank-Wolfe gap at $x$ is
\[
g_{\mathrm{FW}}^{\mathrm{Cl}}(x;\, F,\, \Omega)
  \;:=\; \min_{g \in \partial F(x)}\;
  \max_{u \in \Omega} \langle -g,\, u - x \rangle.
\]
\end{definition}

We are now ready to introduce the concept of Nearby-Clarke Frank-Wolfe Stationarity.

\begin{definition}[Nearby-Clarke Frank-Wolfe Stationary Point (NCFWSP)]
\label{def:ncfwsp} Under the same assumptions on $F$ and $\Omega$, a point
$x \in \Omega$ is a $(\delta, \epsilon)$-NCFWSP if
\[
\exists \tilde{x} \in \mathbb{B}_\delta(x) \cap \Omega \quad \text{s.t.} \quad g_{\mathrm{FW}}^{\mathrm{Cl}}(\tilde{x};\, F,\, \Omega) \le \epsilon.
\]
\end{definition}

When $F$ is continuously differentiable, the Clarke
Frank-Wolfe gap coincides with the classical Frank-Wolfe
gap. The nearby-point definition therefore reduces to
approximate smooth stationarity at a nearby feasible point,
and to the usual condition at $x$ when $\delta=0$.

\begin{remark}
An alternative stationarity definition for the nonsmooth constrained case, namely the Goldstein Frank-Wolfe Stationary Point (GFWSP), can be found in \cite{liu2024zeroth} as a point satisfying
\[
\min_{g \in \partial_\delta F(x)} \; \max_{u \in \Omega} \; \langle -g,\, u - x \rangle ,
\]
where, given $\delta > 0$, the gap is measured by considering the $\delta$-Goldstein subdifferential \cite{goldstein1977optimization}, defined as
\[
\partial_\delta F(x) := \mathrm{conv}\!\left(\bigcup_{y \in \mathbb{B}_\delta(x)} \partial F(y)\right).
\]
In simpler words, in \cite{liu2024zeroth} the gap is measured at the point $x$ itself using an element of the Goldstein subdifferential, which includes subgradients at nearby points. In our opinion, this might be excessively restrictive while the definition we give in this paper is conceptually more suitable for a constrained setting and is a more direct generalization of the Clarke Frank-Wolfe (CFW) gap. In fact, in this paper, for a point to be stationary it suffices that it is in the $\delta$-neighborhood of a point $\tilde{x}$ that is CFW stationary.
\end{remark}

Since we are working with Clarke subdifferentials, the analysis in this subsection relies on a variant of Algorithm~\ref{alg:blds-simplex} in which the finite set of Frank-Wolfe directions $\mathcal{D}^k_\Delta = \{\pm(e_i - w^k)\}_{i=1}^{m}$ is replaced by a \emph{dense} set of directions in the tangent cone to the simplex at $w^k$. In practice, only Lines~3-6 of Algorithm~\ref{alg:blds-simplex} change: instead of cycling over atoms, the algorithm draws directions from a dense subset of $\mathrm{cone}\{e_i - w^k\}_{i=1}^{m}$. Crucially, the maximum feasible stepsize $\bar{\alpha}(w^k, d)$ still admits a closed-form expression for any direction $d$ in the tangent cone to the simplex.

\begin{assumption}[Dense Direction Set]
\label{ass:dense_directions}
During every iteration $k$, the set of directions considered by the algorithm at lines~3-6 is dense in the cone generated by the elementary Frank-Wolfe directions $\{e_i - w^k\}_{i=1}^m$, denoted as $cone(\mathcal{D}^k_\Delta)$.
\end{assumption}

We first establish that, with $\alpha_{\min} > 0$, the iterates become eventually constant.

\begin{lemma}
\label{lem:eventual_constancy}
Under Assumptions~\ref{ass:lls}-\ref{ass:lipschitz_f_x}, with $\gamma > 0$ and $\alpha_{\min} > 0$, the sequence $\{w^k\}$ produced by Algorithm~\ref{alg:blds-simplex} is eventually constant: there exists $\bar{k}$ such that $w^{k+1} = w^k$ for all $k \ge \bar{k}$.
\end{lemma}

\begin{proof}
Every successful iteration decreases $\Ftilde_A$ by at least $\gamma \alpha_{\min}^2$. Since $\Ftilde_A$ is bounded below by $\tilde{f}_{\min}$, the number of successful iterations is finite. After the last successful iteration, all subsequent iterations are unsuccessful and the iterate remains unchanged.
\end{proof}

\begin{theorem}[Nonsmooth convergence to NCFWSP]
\label{thm:gfw_convergence}
Let Assumptions~\ref{ass:lls}--\ref{ass:lipschitz_f_x}
and~\ref{ass:dense_directions} hold, with $\gamma > 0$
and $\alpha_{\min} > 0$. Define
$\eta := 2^{\frac{3}{2}}
\big(\gamma\alpha_{\min}
+2L_f\varepsilon/\alpha_{\min}\big)$.
Then every limit point $w^*$ of the sequence $\{w^k\}$
is an $(\alpha_{\min},\eta)$-Nearby-Clarke Frank-Wolfe
stationary point for $F_A$ on $\Delta^{m-1}$:
there exists
$\tilde{w}^* \in
\mathbb{B}_{\alpha_{\min}}(w^*)\cap\Delta^{m-1}$
such that
\begin{equation}
\label{eq:gfwsp_bound}
\min_{g \in \partial F_A(\tilde{w}^*)} \; \max_{u \in \Delta^{m-1}} \langle -g,\, u - \tilde{w}^* \rangle \le 2^{\frac{3}{2}} \Big(\gamma \alpha_{\min} + \frac{2L_f \varepsilon}{\alpha_{\min}}\Big).
\end{equation}
\end{theorem}

\begin{proof}
By Lemma~\ref{lem:eventual_constancy} the sequence \(\{w^k\}\) is eventually constant and we denote its limit value with \(w^*\). For \(k\) large enough, every direction trial is unsuccessful and \(\alpha^k=\alpha_{\min}\), therefore for every \(d \in cone(\mathcal{D}^*_{\Delta})\) the effective stepsize used by the algorithm is:
\[
\hat{\alpha}^*(d) := \min\{\alpha_{\min}, \bar{\alpha}(w^*,d)\}.
\]

By the unsuccessful condition, we have:

\[
\tilde F_A(w^*+ \hat{\alpha}^*(d) d) \;>\; \tilde F_A(w^*) - \gamma\alpha_{\min}^2.
\]

Using Remark~\ref{rem:bounded_objective_error} (bounded approximation error \(|\tilde F_A - F_A|\le L_f\varepsilon\)) we obtain, for every \(d \in cone(\mathcal{D}^*_{\Delta})\),
\begin{equation}
    \label{eq:unsuccessful_asymptotic}
F_A(w^*) \;<\; F_A(w^*+ \hat{\alpha}^*(d) d) + \gamma\alpha_{\min}^2 + 2L_f\varepsilon.
\end{equation}
We call the tangent cone to the simplex at a point $w$ as $T_{\Delta^{m-1}}(w)$. 
Define the set
$$
\mathcal{H}
:= \{h=\alpha_{\min}d:
d\in T_{\Delta^{m-1}}(w^*),\ \|d\|=1\}
$$
of tangent displacements with norm $\alpha_{\min}$, and the set
$$
\mathcal{G}
:= \{g=\hat{\alpha}^*(d)d:
d\in T_{\Delta^{m-1}}(w^*),\ \|d\|=1\}
$$
of clipped displacements for which $w^*+g$ is feasible.
Let $\chi_{\Omega}(x)$ be the characteristic function of the set $\Omega$, i.e. $\chi_{\Omega}(x) = 0$ if $x \in \Omega$ and $+\infty$ otherwise.
From \eqref{eq:unsuccessful_asymptotic} and density of the generated directions it holds that, for every $h \in \mathcal{H}$,
    \[
    F_A(w^*) \;<\; F_A(w^*+h) \;+\; C \;+\; \chi_{\Delta^{\bar{m}-1}}(w^*+h),
    \qquad C:=\gamma\alpha_{\min}^2+2L_f\varepsilon,
    \]
    where the penalty $\chi_{\Delta^{\bar{m}-1}}(w^*+h) = +\infty$ is given when $w^* + h$ is not feasible.

    Let us define the auxiliary function:
    \[
    \bar{F}_x(h) \;:=\; F_A(x+h) \;+\; C \frac{\|h\|^2}{\alpha_{\min}^2} \;+\; \chi_{\Delta^{\bar{m}-1}}(w^*+h),
    \]

    and note that 
    \begin{equation}
        \label{eq:auxiliary_suboptimality}
    \bar{F}_{w^*}(0) < \bar{F}_{w^*}(h) \quad \text{for every } h \in \mathcal{H}.
    \end{equation}

    We now consider the compact sets bounded by $\mathcal{H}$ and $\mathcal{G}$ and centered at the origin, namely: 
    \begin{align*}
    \bar{\mathcal{H}} &:= \{\tau h \mid h \in \mathcal{H}, \tau \in [0,1]\} \subseteq \mathbb{B}_{\alpha_{\min}}(0), \\
    \bar{\mathcal{G}} &:= \{\tau h \mid h \in \mathcal{G}, \tau \in [0,1]\} \subseteq \bar{\mathcal{H}}.
    \end{align*}
    Figure~\ref{fig:diagramHG} provides a diagram which can help understand the difference between the four different sets involved here. Notice that $\bar{F}_x(h) = +\infty$ for every $h \in \bar{\mathcal{H}} \backslash \bar{\mathcal{G}}$, and that $\{w^* + h \mid h \in \bar{\mathcal{G}}\} \subseteq \Delta^{\bar m-1}$.

    By compactness of $\bar{\mathcal{G}}$ and the continuity of $F_A$, there exists a minimizer $\tilde{h} \in \bar{\mathcal{G}}$ of $\bar{F}_{w^*}$ which is also a minimizer inside $\bar{\mathcal{H}}$. Since the boundary of $\bar{\mathcal{H}}$ is suboptimal from \eqref{eq:auxiliary_suboptimality}, thus $\|\tilde{h}\|<\alpha_{\min}$, i.e. the minimizer is an interior point of $\bar{\mathcal{H}}$ and $w^* + \tilde{h} \in \Delta^{\bar m - 1}$ because $\tilde{h} \in \bar{\mathcal{G}}$. 
    It holds the following Clarke Frank-Wolfe optimality condition:
    \begin{equation}
        \label{eq:clarke_fw_optimality}
        \min_{g \in \partial F_{w^*}(\tilde{h})} \; \max_{d \in T_{\bar{\mathcal{G}}}(w^* + \tilde{h})} \langle -g, d \rangle = 0.
    \end{equation}

    Moreover, it holds that:
    \begin{equation}
        \label{eq:tangent_cones_equality}
    T_{\bar{\mathcal{G}}}(\tilde{h}) = T_{\Delta^{\bar m -1}}(w^* + \tilde{h})
    \end{equation}
    because $\{w^* + \mathcal{G}\} \subset \Delta^{\bar{m}-1}$ and if $\tilde{h}$ is in the border of $\mathcal{G}$, then $w^* + \tilde{h}$ is in the border of $\Delta^{\bar{m} -1}$.

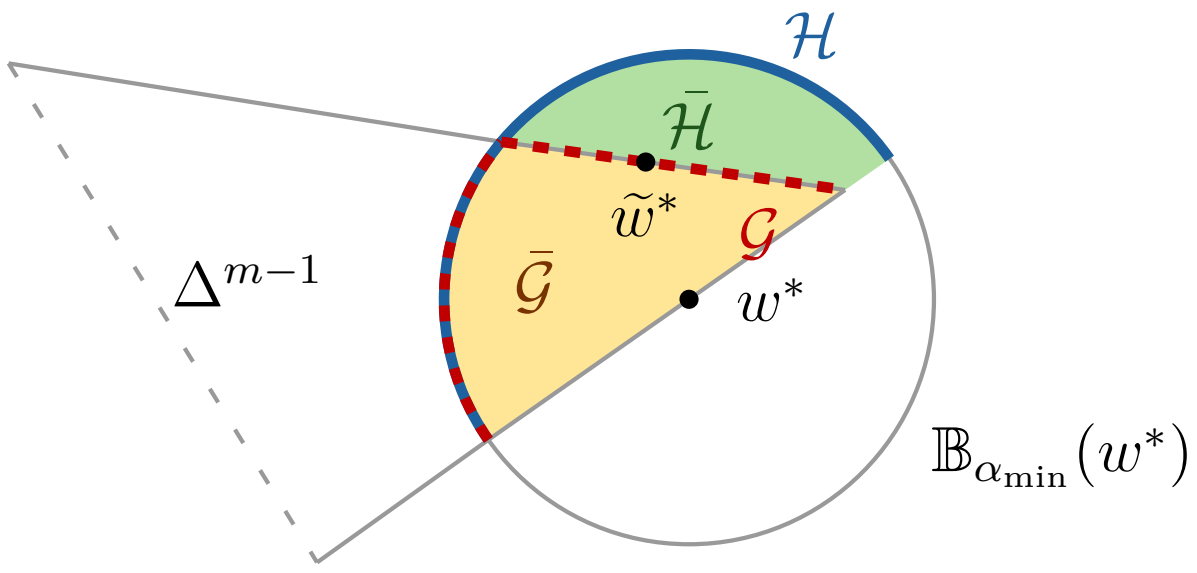
\begin{figure}
\centering
\begin{tikzpicture}[scale=0.6]

    \def\R{2.7}
    \def\angBot{215}
    \def\angTop{140}
    \def\angDiag{35}
    \def\rApex{2.1}

    \coordinate (Wstar) at (0,0);
    \coordinate (Pbot) at (\angBot:\R);
    \coordinate (Ptop) at (\angTop:\R);
    \coordinate (Pblue) at (\angDiag:\R);
    \coordinate (Apex) at (\angDiag:\rApex);

    \coordinate (Vtop) at (-7.5,2.6);
    \coordinate (Vbot) at (-4.1,-2.9);

    \fill[mygreenfill]
        (Pbot) arc (\angBot:\angDiag:\R) -- cycle;

    \fill[myorangefill]
        (Pbot) arc (\angBot:\angTop:\R)
        -- (Apex) -- cycle;

    \draw[gray!80,line width=0.8pt]
        (Vtop) -- (Ptop) -- (Apex);
    \draw[gray!80,line width=0.8pt]
        (Vbot) -- (Pbot) -- (Apex);
    \draw[gray!80,line width=0.8pt,loosely dashed]
        (Vtop) -- (Vbot);

    \draw[gray!80,line width=0.8pt]
        (Wstar) circle (\R);

    \draw[myblue,line width=2.0pt]
        (Pbot) arc (\angBot:\angDiag:\R);

    \draw[myred,line width=2.0pt,dashed]
        (Pbot) arc (\angBot:\angTop:\R);
    \draw[myred,line width=2.0pt,dashed]
        (Ptop) -- (Apex);

    \fill[black] (Wstar) circle (3pt);

    \coordinate (Wtilde) at ($(Ptop)!0.42!(Apex)$);
    \fill[black] (Wtilde) circle (3pt);

    \node[left] at (-3.8,0.2)
        {\large $\Delta^{m-1}$};

    \node[goldenscript] at (-1.7,0.2)
        {\large $\bar{\mathcal G}$};

    \node[myred,anchor=north]
        at ($(Ptop)!0.75!(Apex)+(0,-0.12)$)
        {\large $\mathcal G$};

    \node[below=2pt] at (Wtilde)
        {\large $\widetilde w^*$};

    \node[right=5pt,anchor=west] at (Wstar)
        {\large $w^*$};

    \node[mygreentext,below=3pt] at (90:\R)
        {\large $\bar{\mathcal H}$};

    \node[myblue,above=2pt] at (60:\R)
        {\large $\mathcal H$};

    \node[anchor=north west] at (-25:\R)
        {\large $\mathbb B_{\alpha_{\min}}(w^*)$};

\end{tikzpicture}
\caption{Example of sets $\mathcal H$ (blue),
$\mathcal G$ (red dashed), $\bar{\mathcal H}$
(green, including the golden region), and
$\bar{\mathcal G}$ (golden) from the proof of
Theorem~\ref{thm:gfw_convergence}.
With a slight abuse of notation, all four sets are shown
translated by the vector $w^*$. Both shaded regions include their boundaries.
The golden region is the intersection of the closed
tangent half-disk with the simplex; in particular,
its straight boundary segments are included in
$\bar{\mathcal G}$, even where they do not belong
to $\mathcal G$.}
\label{fig:diagramHG}
\end{figure}
    Therefore, from \eqref{eq:clarke_fw_optimality} and \eqref{eq:tangent_cones_equality}, and the fact that $\partial F_{w^*}(\tilde{h}) = \partial F_A(w^* + \tilde{h}) + \tfrac{2C}{\alpha_{\min}^2} \tilde{h}$, it follows that:
    \begin{align}
        & \min_{g \in \partial F_A(w^* + \tilde{h})} \; \max_{u \in \Delta^{\bar{m}-1}} \langle -g - \tfrac{2C}{\alpha_{\min}^2} \tilde{h}, u - (w^*+\tilde{h}) \rangle = 0 \\
        \label{eq:intermediate_gfw_bound}
        \Rightarrow
        & \min_{g \in \partial F_A(w^* + \tilde{h})} \; \max_{u \in \Delta^{\bar{m}-1}} \left( \langle -g, u - (w^*+\tilde{h}) \rangle - \frac{2C}{\alpha_{\min}^2} \langle \tilde{h}, u - (w^*+\tilde{h}) \rangle \right) = 0.
    \end{align}

    Now, we know that for general functions $A(x)$ and $B(x)$, we have $\max_x (A(x) - B(x)) = 0 \Rightarrow A(x) \leq B(x) \quad \forall x$; in particular $\max_x A(x) \leq B(\text{arg}\max_x(A(x))) \leq \max_x B(x)$. Using this property, from \eqref{eq:intermediate_gfw_bound} we imply that there exists $\bar{g} \in \partial F_A(w^* + \tilde{h})$ such that:

    \begin{align*}
        & \max_{u \in \Delta^{\bar{m}-1}} \left( \langle -\bar{g}, u - (w^*+\tilde{h}) \rangle - \frac{2C}{\alpha_{\min}^2} \langle \tilde{h}, u - (w^*+\tilde{h}) \rangle \right) = 0 \\
        \Rightarrow
        & \max_{u \in \Delta^{\bar{m}-1}} \langle -\bar{g}, u - (w^*+\tilde{h}) \rangle \leq \max_{u \in \Delta^{\bar{m}-1}} \frac{2C}{\alpha_{\min}^2} \langle \tilde{h}, u - (w^*+\tilde{h}) \rangle \leq \frac{2C}{\alpha_{\min}^2} \|\tilde{h}\| D.
    \end{align*}

    Finally, using $D=\sqrt{2}$, it follows that:
    \begin{align*}
        & \min_{g \in \partial F_A(w^* + \tilde{h})} \; \max_{u \in \Delta^{\bar{m}-1}} \langle -g, u - (w^* + \tilde{h}) \rangle \leq \frac{2C}{\alpha_{\min}^2} \|\tilde{h}\| D \\
        \Rightarrow
        & \;\;\exists \tilde{w}^* = w^* + \tilde{h} \in \mathbb{B}_{\alpha_{\min}}(w^*) \text{ s.t. } \min_{g \in \partial F_A(\bar{w^*})} \; \max_{u \in \Delta^{m-1}} \langle -g, u - \tilde{w}^* \rangle \leq \frac{2^{\frac{3}{2}} C}{\alpha_{\min}}.
    \end{align*}
\end{proof}

\section{Numerical Experiments}
\label{sec:numerical_experiments}

In this section we present our results validating BLDS-SIMPLEX (Algorithm~\ref{alg:blds-simplex}) on bilevel benchmark problems from the BOLIB library~\cite{zhou2020bolib,bolib3}, showing that BLDS-SIMPLEX is always able to find high-quality solutions. As for number of function evaluations, we note that in worst case scenarios it must check $O(m)$ directions at each iteration. To make the algorithm more scalable and competitive, we also analyze the performance of a variant in which we take a random subset of FW directions at each iteration, instead of the full batch. The variant searches randomly ordered blocks of Frank--Wolfe directions and stops the search as soon
as a sufficient-decrease step is found. This subspace variant BLDS-SIMPLEX-S, described in Algorithm~\ref{alg:blds-simplex-subspaces}, reduces in fact the computational cost per iteration without sacrificing solution quality, making it a promising choice in practice. The results are compared with well known solvers for both bilevel and single-level derivative-free optimization. Experiments can be reproduced with the code available in \href{https://github.com/GibboSan/df-blo}{GitHub} \footnote{https://github.com/GibboSan/df-blo}.

\paragraph{Test problems and polytope constraint construction.}
We select nine bilevel problems with dimensions $n_x$ ranging from $2$ to $10$ and $n_y$ from $2$ to $12$, covering linear, quadratic, and nonlinear objectives; see Table~\ref{tab:problems}.
For each problem, the upper-level polytope $\Omega_x = \mathrm{conv}(a_1,\dots,a_m)$ is constructed by sampling $m$ vertices uniformly on the surface of a ball of radius $r=5$ in $\reals^{n_x}$. 
The polytope is then shifted so that the known optimal solution $x^*$ is positioned at the centroid.

\begin{table}[t]
\centering
\caption{Test problems from BOLIB~\cite{bolib3,zhou2020bolib} Objective types at upper level (UL) and lower level (LL): L\,=\,linear, Q\,=\,quadratic, NL\,=\,nonlinear.}
\label{tab:problems}
\begin{tabular}{l c c c c}
\hline
\textbf{Problem} & $n_x$ & $n_y$ & \textbf{UL} & \textbf{LL} \\
\hline
DempeLohse2011Ex31a   & 2  & 2  & Q  & L  \\
CalamaiVicente1994b   & 4  & 2  & Q  & Q  \\
CalamaiVicente1994c   & 4  & 2  & Q  & Q  \\
SinhaMaloDeb2014TP9   & 10 & 10 & Q  & NL \\
SinhaMaloDeb2014TP10  & 10 & 10 & Q  & NL \\
NieWangYe2017Ex54     & 4  & 4  & NL & Q  \\
hypercube\_10         & 10 & 10 & L  & L  \\
qbp\_8\_12            & 8  & 12 & Q  & Q  \\
qbp\_9\_9             & 9  & 9  & Q  & Q  \\
\hline
\end{tabular}
\end{table}

Each of the $9$ problems is combined with $m\in\{50,100,200,500,1000\}$ atoms and $3$ random seeds ($\{11,22,33\}$), for a total of $9\times5\times3=135$ instances per solver ($27$ per fixed $m$).
In the discussion we group vertex counts into \emph{small} instances ($m\in\{50,100,200\}$) and \emph{large} instances ($m\in\{500,1000\}$) to highlight scalability trends.

\paragraph{Lower-level oracle.}
For each upper-level trial point $x$, the lower-level problem is solved by SciPy's SLSQP \cite{kraft1988software,2020SciPy-NMeth} with a fixed tolerance of $10^{-3}$. We stress that no rigorous certificate of attaining a prescribed accuracy $\varepsilon$ is available at every call; SLSQP is used as a practical black-box oracle with a fixed stopping criterion.

\paragraph{Subspace variant.}
Cycling through all $2m$ Frank-Wolfe directions per iteration becomes expensive. Algorithm~\ref{alg:blds-simplex-subspaces} mitigates this by randomly partitioning $\{1,\dots,m\}$ into blocks of size $s$ and searching one block at a time: upon finding a sufficient-decrease step the iteration ends immediately; otherwise the stepsize is contracted and the next block is tried. We set $s = n_x$, keeping the per-block cost proportional to the original problem dimension. When $s=m$ the algorithm reduces to the original. The same convergence  guarantees hold. In the experiment we will report the subspace variant as \textbf{BLDS} and the original strategy as \textbf{BLDS (full)}.

\begin{algorithm}[t]
\caption{BLDS-SIMPLEX-S}
\label{alg:blds-simplex-subspaces}
\begin{algorithmic}[1]
\Require $w^0\!\in\!\Delta^{m-1}$, $\alpha^0$, $\alpha_{\min}$, $\alpha_{\max}$, $\theta,\delta\!\in\!(0,1)$, $\gamma\!>\!0$, block size $s\!\geq\!1$
\For{$k = 0,1,2,\ldots$}
    \State $\texttt{success}\leftarrow\texttt{false}$;\; $\alpha \leftarrow \max\{\alpha^k,\alpha_{\min}\}$
    \State Randomly partition $\{1,\dots,m\}$ into blocks $\mathcal{B}_1,\dots,\mathcal{B}_{\lceil m/s\rceil}$
    \For{$b = 1,\dots,\lceil m/s\rceil$}
        \For{each $i\in\mathcal{B}_b$ and $\mathrm{sgn}\in\{+,-\}$}
            \State $d\leftarrow\mathrm{sgn}\,(e_i - w^k)$;\; $\hat{\alpha} \leftarrow\min\{\alpha,\bar\alpha(w^k,d)\}$
            \If{$\hat{\alpha}>0$ \textbf{and} $\Ftilde_A(w^k\!+\!\hat{\alpha} d)\le\Ftilde_A(w^k)-\gamma(\alpha)^2$}
                \State $w^{k+1}\leftarrow w^k+\hat{\alpha} d$;\; $\texttt{success}\leftarrow\texttt{true}$;\; \textbf{break}
            \EndIf
        \EndFor
        \If{$\texttt{success}$}
            $\alpha \leftarrow \min\{\alpha_{\max},\,\alpha/\delta\}$;\; \textbf{break}
        \Else\; $\alpha \leftarrow \max\{\alpha_{\min},\,\theta\alpha\}$
        \EndIf
    \EndFor
    \State $\alpha^{k+1} \leftarrow \alpha$.
    \If{not $\texttt{success}$ and $\alpha^{k+1} =\alpha^k$} \textbf{Terminate.}
    \EndIf
\EndFor
\end{algorithmic}
\end{algorithm}

\paragraph{Competing algorithms and parameters.} We compare BLDS-SIMPLEX and BLDS-SIMPLEX-S against five derivative-free solvers, all using the simplex formulation in \eqref{eq:blop_simplex_reduced} and the same lower-level oracle:
\textbf{DFN-LLA}~\cite{cesaroni2026derivativefreebileveloptimizationinexact} (based on an exact penalty function),
\textbf{DSPFD}~\cite{gratton2019direct} (using probabilistic feasible directions, with the original MATLAB implementation available \href{https://www.lamsade.dauphine.fr/~croyer/numerics.html}{here}),
\textbf{LINCOA}~\cite{powell2015fast} (a model-based trust-region method from the PDFO library \cite{ragonneau2021pdfo}),
\textbf{NOMAD}~\cite{AuDe2006,audet2022Algo} (a mesh adaptive direct-search method from \href{https://pypi.org/project/PyNomadBBO/}{PyNomadBBO}), and
\textbf{ZOFW}~\cite{sahu2019towards} (zeroth-order Frank-Wolfe with I-RDSA gradient estimation).

The single-level solvers DSPFD, LINCOA, NOMAD, and ZOFW are applied directly to the inexact reduced objective
$\Ftilde_A(w)=\Ftilde(Aw)=f(Aw,\tilde y(Aw))$
over the simplex. These solvers do not explicitly exploit the bilevel structure: they treat the values returned by the common lower-level oracle as black-box objective evaluations. Thus, they operate as single-level methods optimizing an inexact objective, whose evaluation errors arise from the approximate solution of the lower-level problem.

All solvers share a budget of $N_{\max}=5000$ inexact objective evaluations.
For our algorithms, we set $\alpha^0=0.4$, $\alpha_{\min}=10^{-5}$, $\alpha_{\max}=0.8$, $\theta=\delta=0.5$, and $\gamma=10^{-3}$.
The initial point is a randomly generated interior point of the polytope, with a fixed seed for reproducibility.

\paragraph{Metrics.}
Performance is assessed via \emph{data profiles}: the fraction of problems solved as a function of the number of function evaluations, where a problem is solved when $f^k \le f^* + \epsilon_{\mathrm{tol}}$. We report results for $\epsilon_{\mathrm{tol}}\in\{10^{-1},\,10^{-3},\,10^{-5}\}$.

\begin{figure}[htbp]
    \centering
    \begin{subfigure}[t]{0.32\linewidth}
        \centering
        \includegraphics[width=\linewidth]{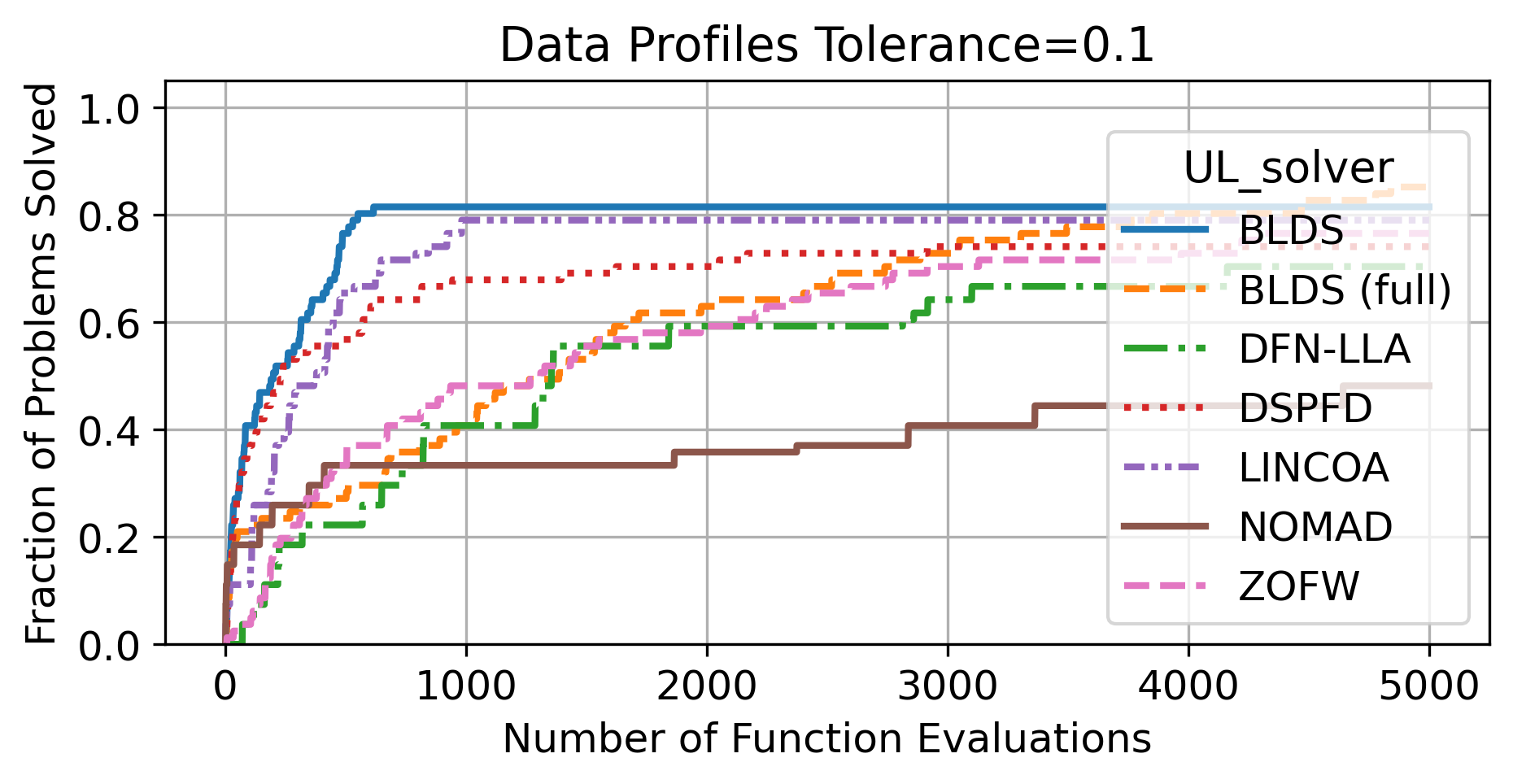}
        \caption{$\epsilon_{\mathrm{tol}}=10^{-1}$}
    \end{subfigure}
    \hfill
    \begin{subfigure}[t]{0.32\linewidth}
        \centering
        \includegraphics[width=\linewidth]{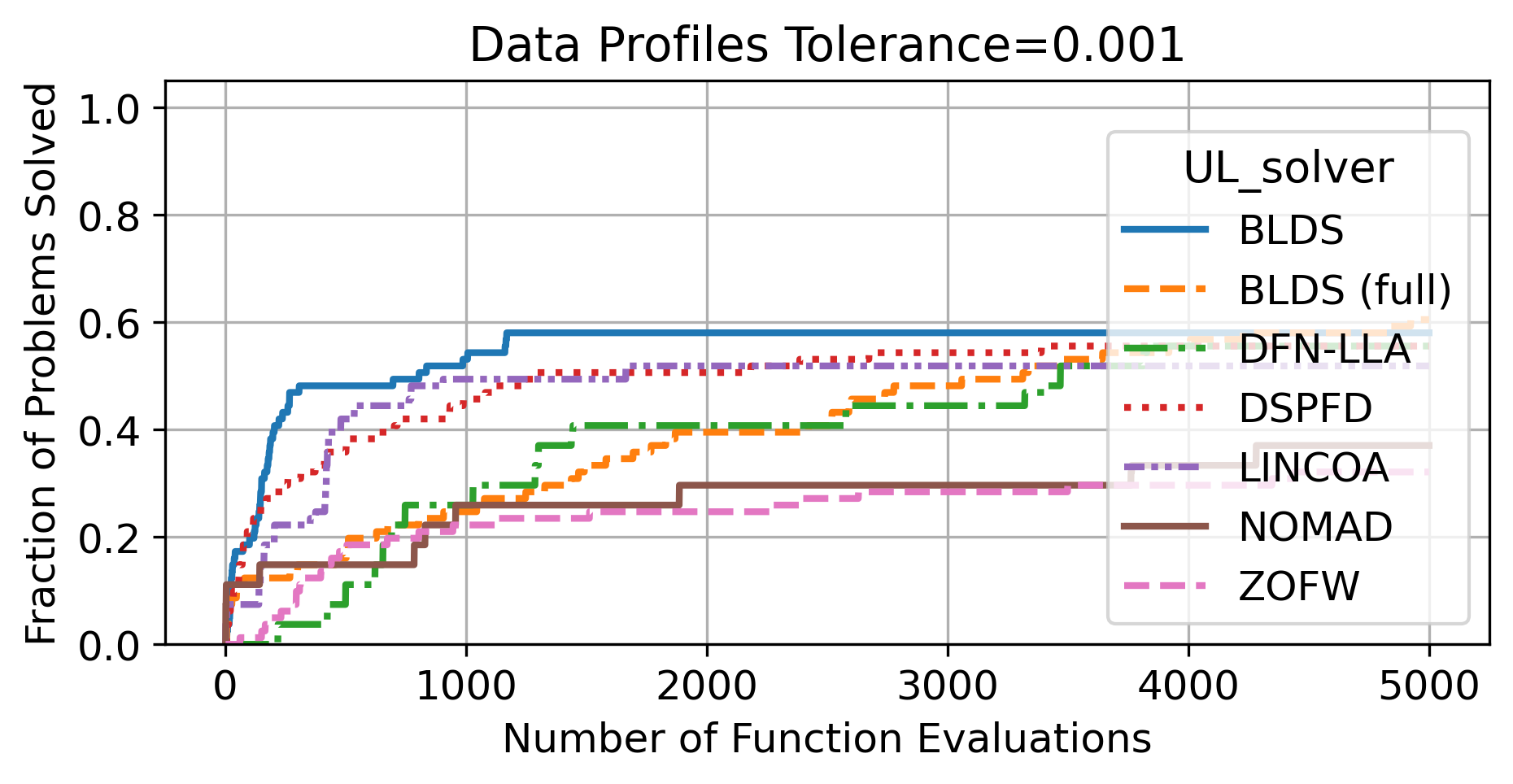}
        \caption{$\epsilon_{\mathrm{tol}}=10^{-3}$}
    \end{subfigure}
    \hfill
    \begin{subfigure}[t]{0.32\linewidth}
        \centering
        \includegraphics[width=\linewidth]{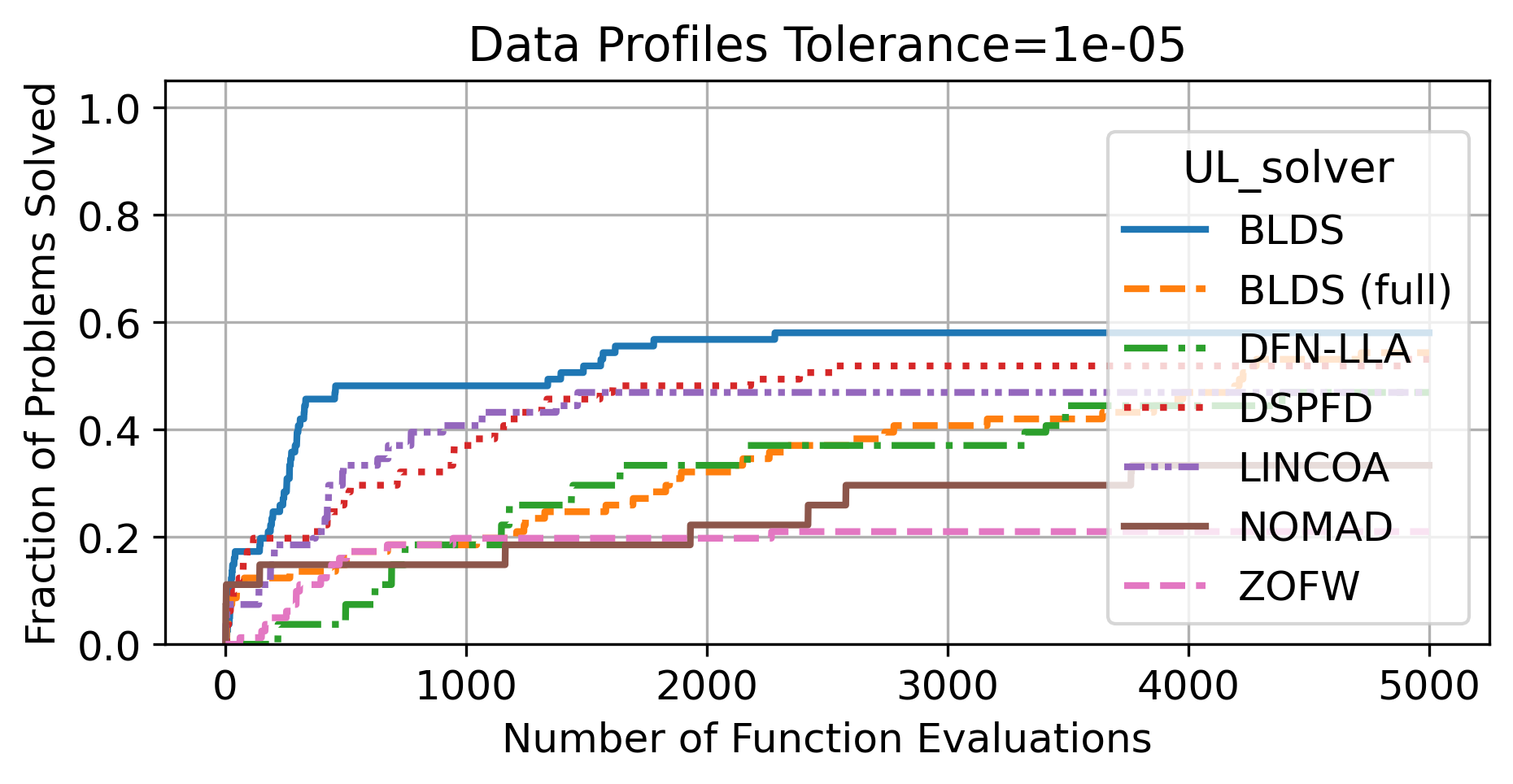}
        \caption{$\epsilon_{\mathrm{tol}}=10^{-5}$}
    \end{subfigure}
    \caption{Data profiles on the whole benchmark set.}
    \label{fig:profiles_all}
\end{figure}

\paragraph{Results.}
We first report data profiles on the whole benchmark set (all values of $m$ aggregated), at three tolerances; see Figure~\ref{fig:profiles_all}. 

From the plots it is clear how the proposed BLDS (subspace variant) is able to outperform all the competitors in terms of efficiency and robustness. The plots also highlight the advantage of the subspace variant over the ``full'' version, whose search over the complete set of FW directions can increase both the per-iteration cost and the total number of function evaluations. Based on these results, we omit ``BLDS (full)'' from the subsequent plots.

\begin{figure}[htbp]
    \centering
    \begin{subfigure}[t]{0.32\linewidth}
        \centering
        \includegraphics[width=\linewidth]{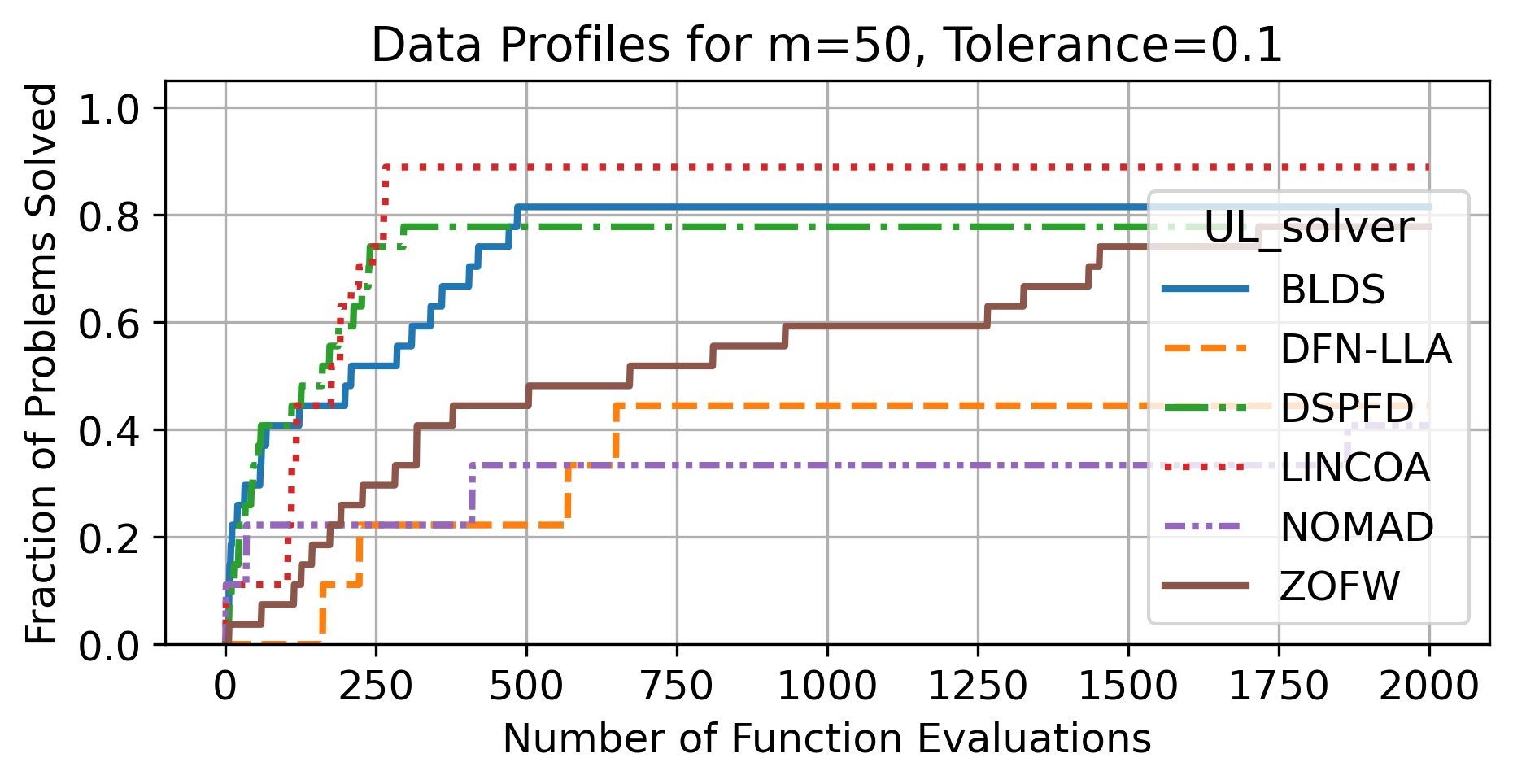}
        \caption{$m=50$}
    \end{subfigure}
    \hfill
    \begin{subfigure}[t]{0.32\linewidth}
        \centering
        \includegraphics[width=\linewidth]{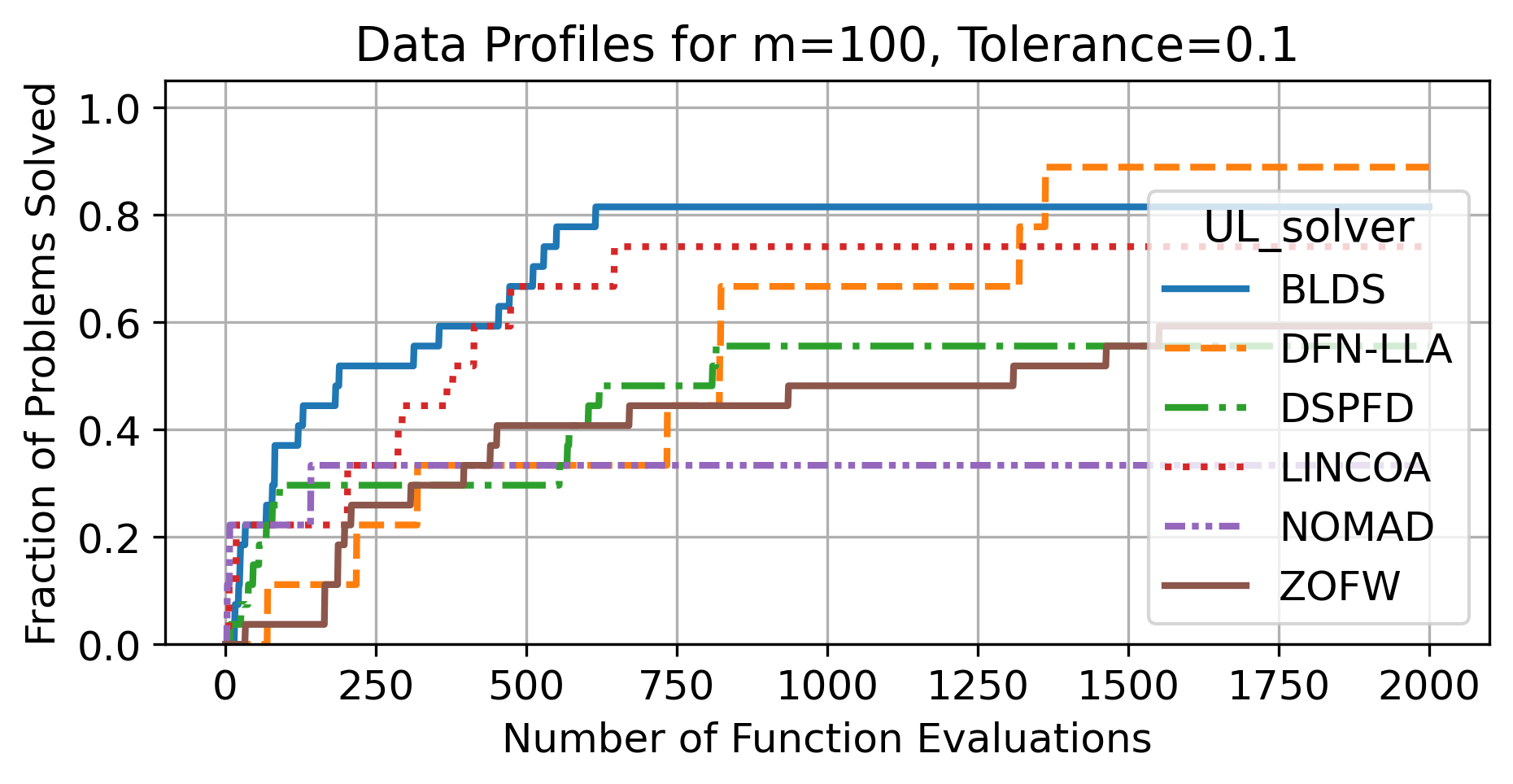}
        \caption{$m=100$}
    \end{subfigure}
    \hfill
    \begin{subfigure}[t]{0.32\linewidth}
        \centering
        \includegraphics[width=\linewidth]{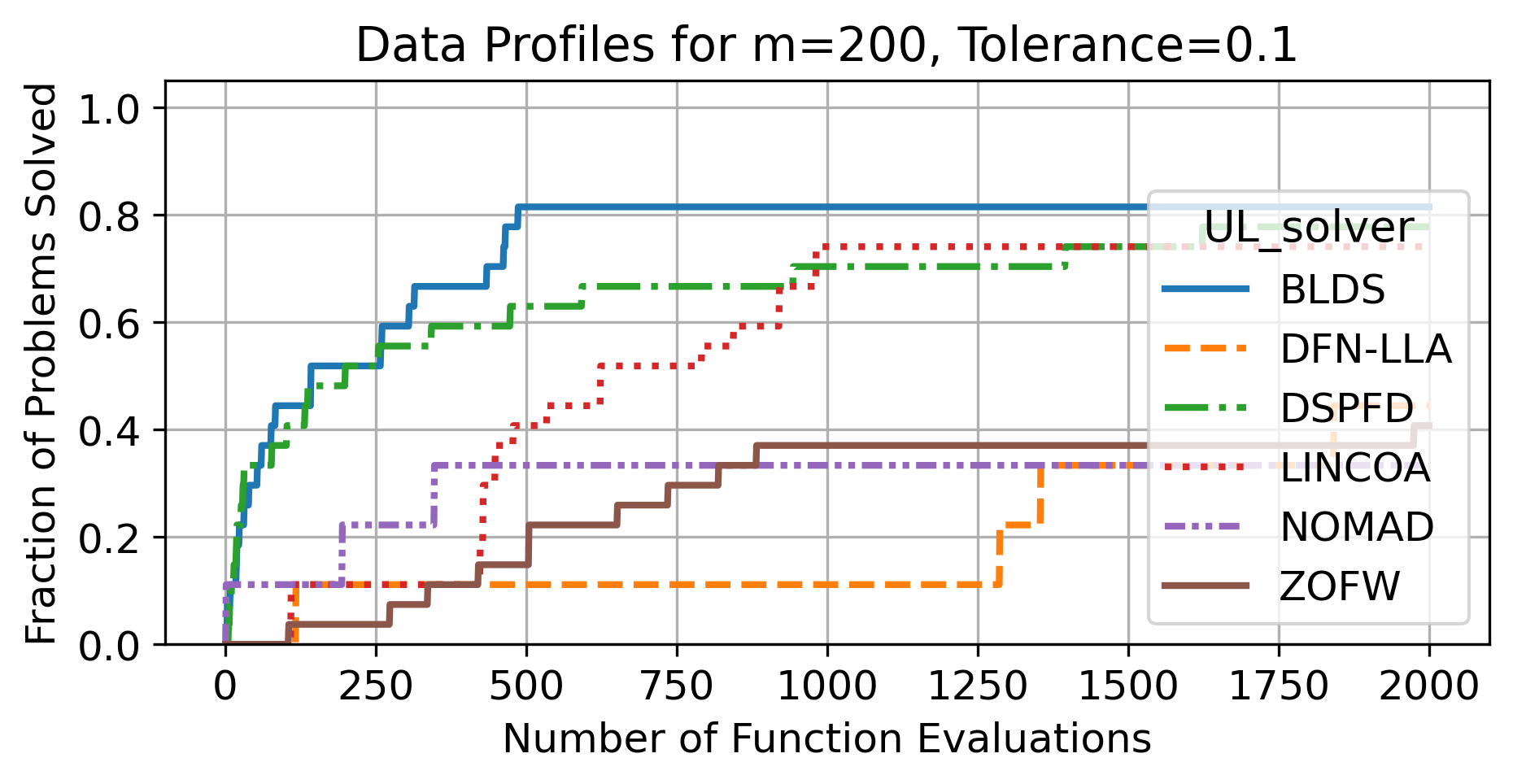}
        \caption{$m=200$}
    \end{subfigure}
    \caption{Data profiles at $\epsilon_{\mathrm{tol}}=10^{-1}$ for small vertex counts.}
    \label{fig:profiles_small_vertices_tol01}
\end{figure}

\begin{figure}[htbp]
    \centering
    \begin{subfigure}[t]{0.32\linewidth}
        \centering
        \includegraphics[width=\linewidth]{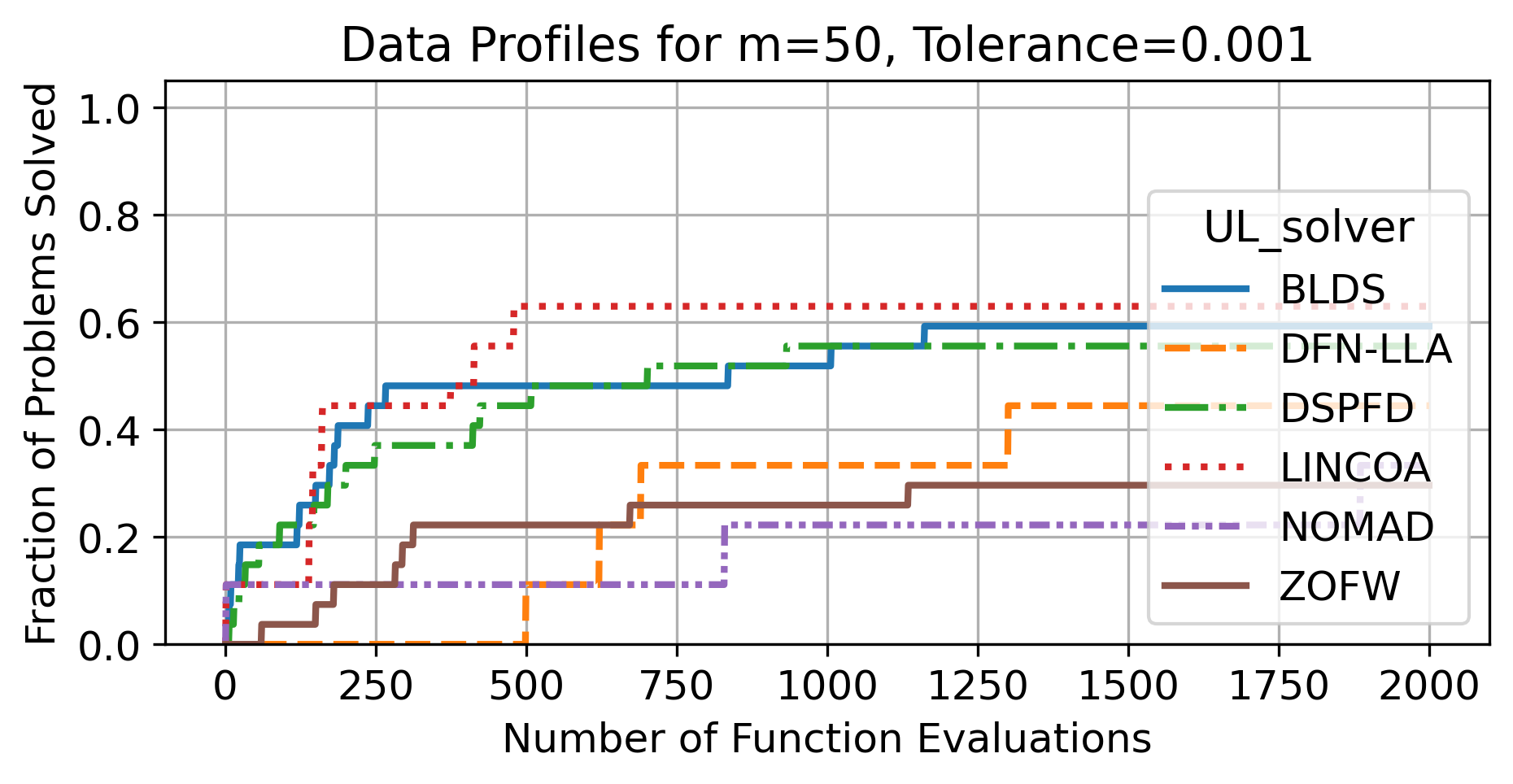}
        \caption{$m=50$}
    \end{subfigure}
    \hfill
    \begin{subfigure}[t]{0.32\linewidth}
        \centering
        \includegraphics[width=\linewidth]{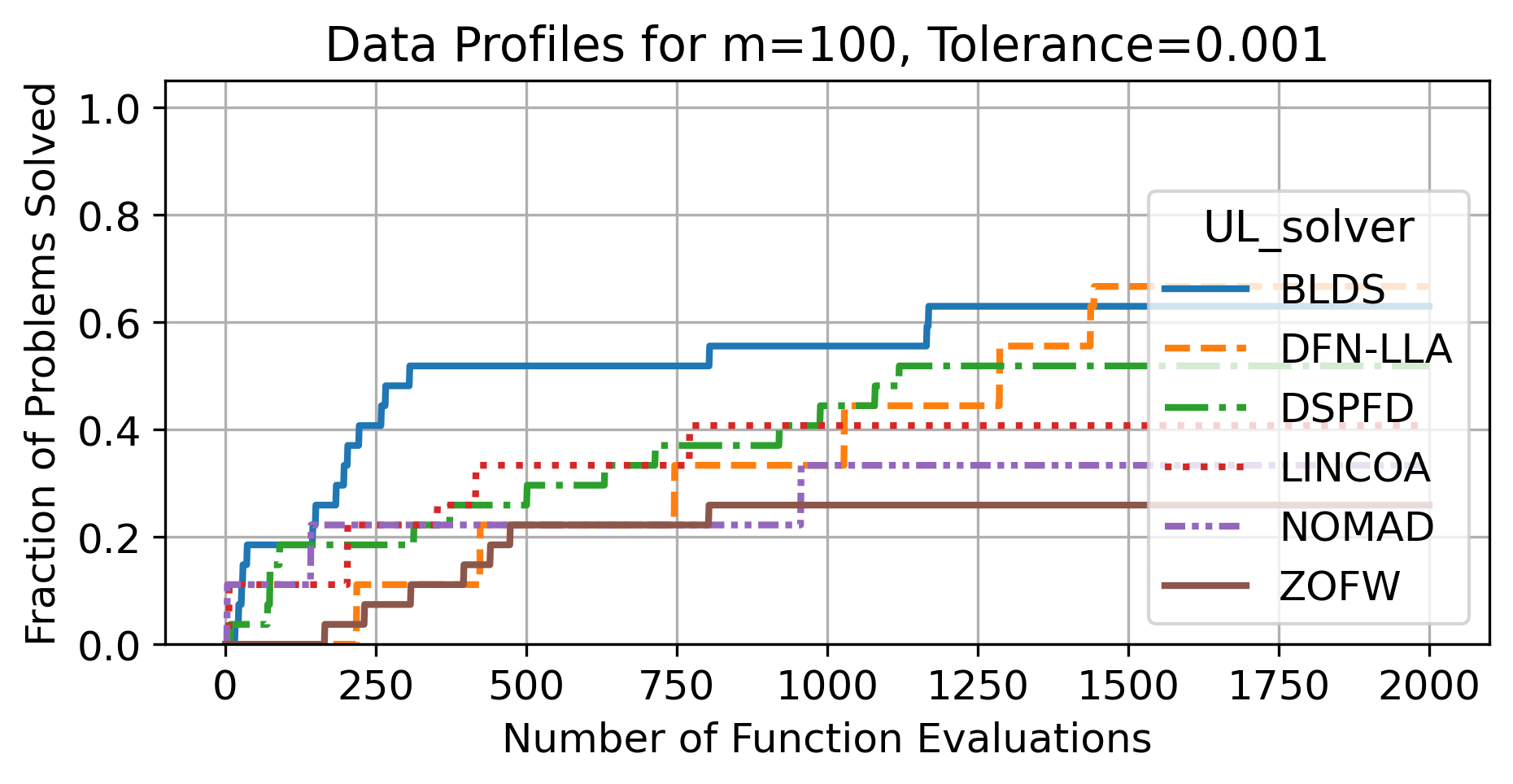}
        \caption{$m=100$}
    \end{subfigure}
    \hfill
    \begin{subfigure}[t]{0.32\linewidth}
        \centering
        \includegraphics[width=\linewidth]{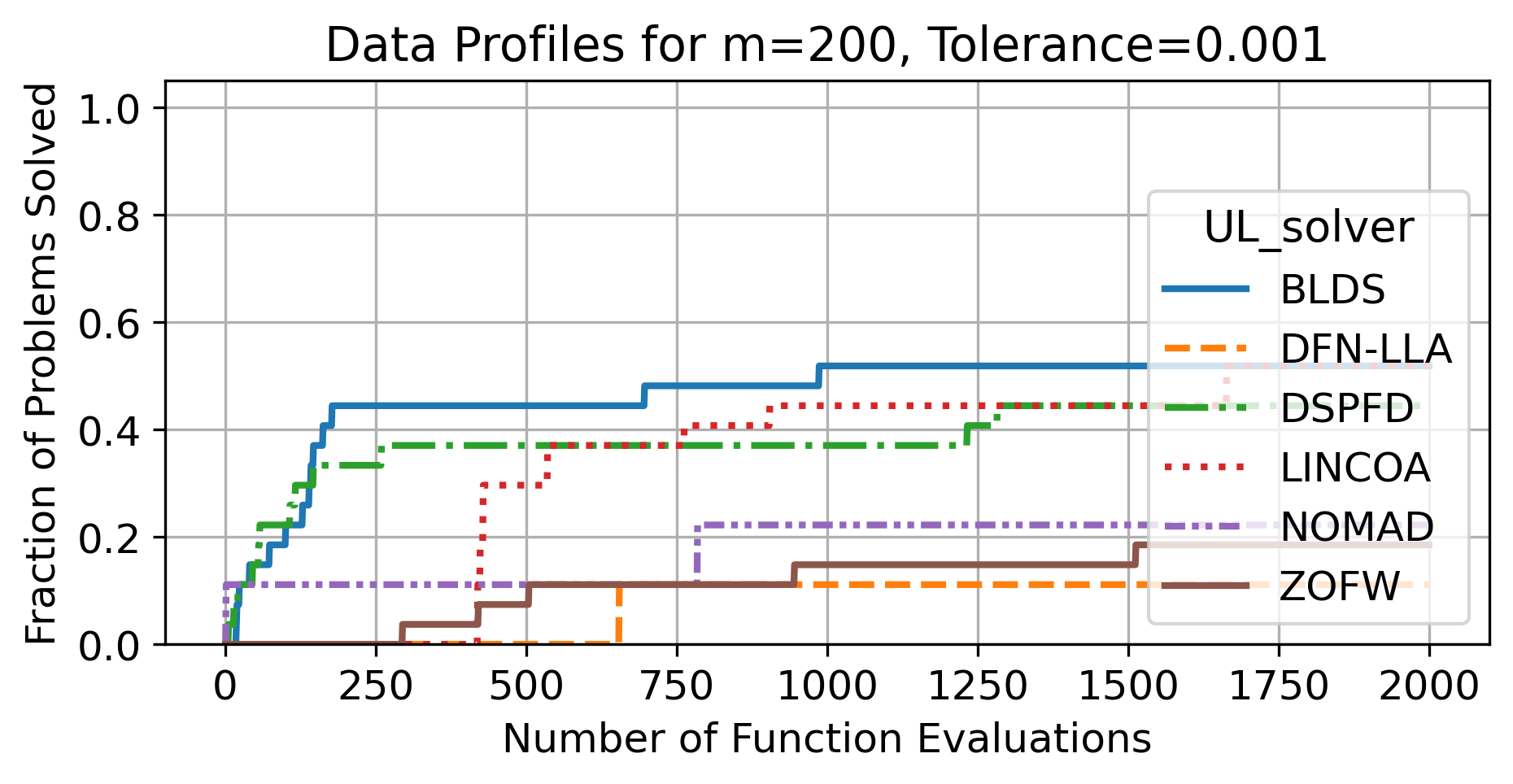}
        \caption{$m=200$}
    \end{subfigure}
    \caption{Data profiles at $\epsilon_{\mathrm{tol}}=10^{-3}$ for small vertex counts.}
    \label{fig:profiles_small_vertices_tol001}
\end{figure}

\begin{figure}[htbp]
    \centering
    \begin{subfigure}[t]{0.32\linewidth}
        \centering
        \includegraphics[width=\linewidth]{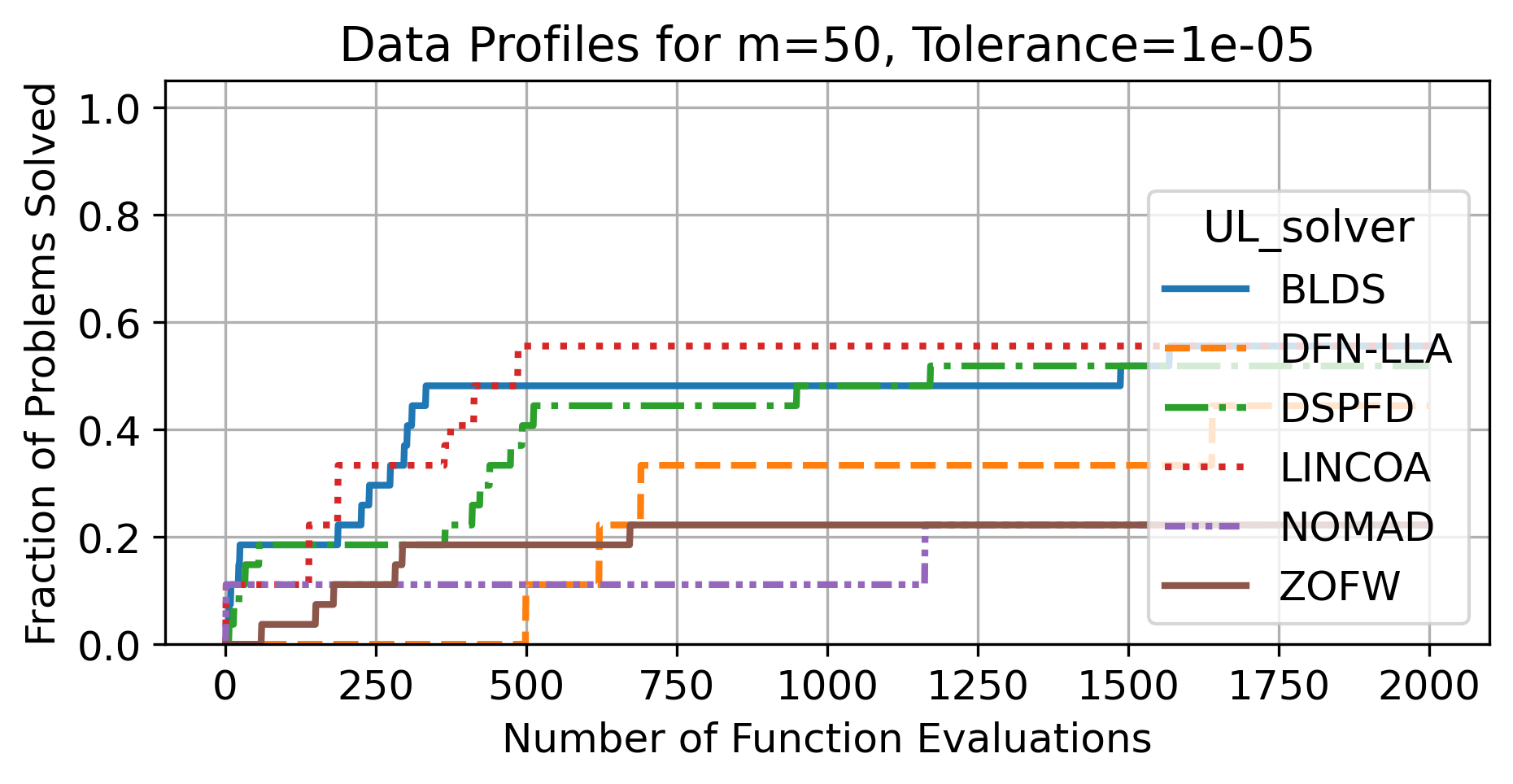}
        \caption{$m=50$}
    \end{subfigure}
    \hfill
    \begin{subfigure}[t]{0.32\linewidth}
        \centering
        \includegraphics[width=\linewidth]{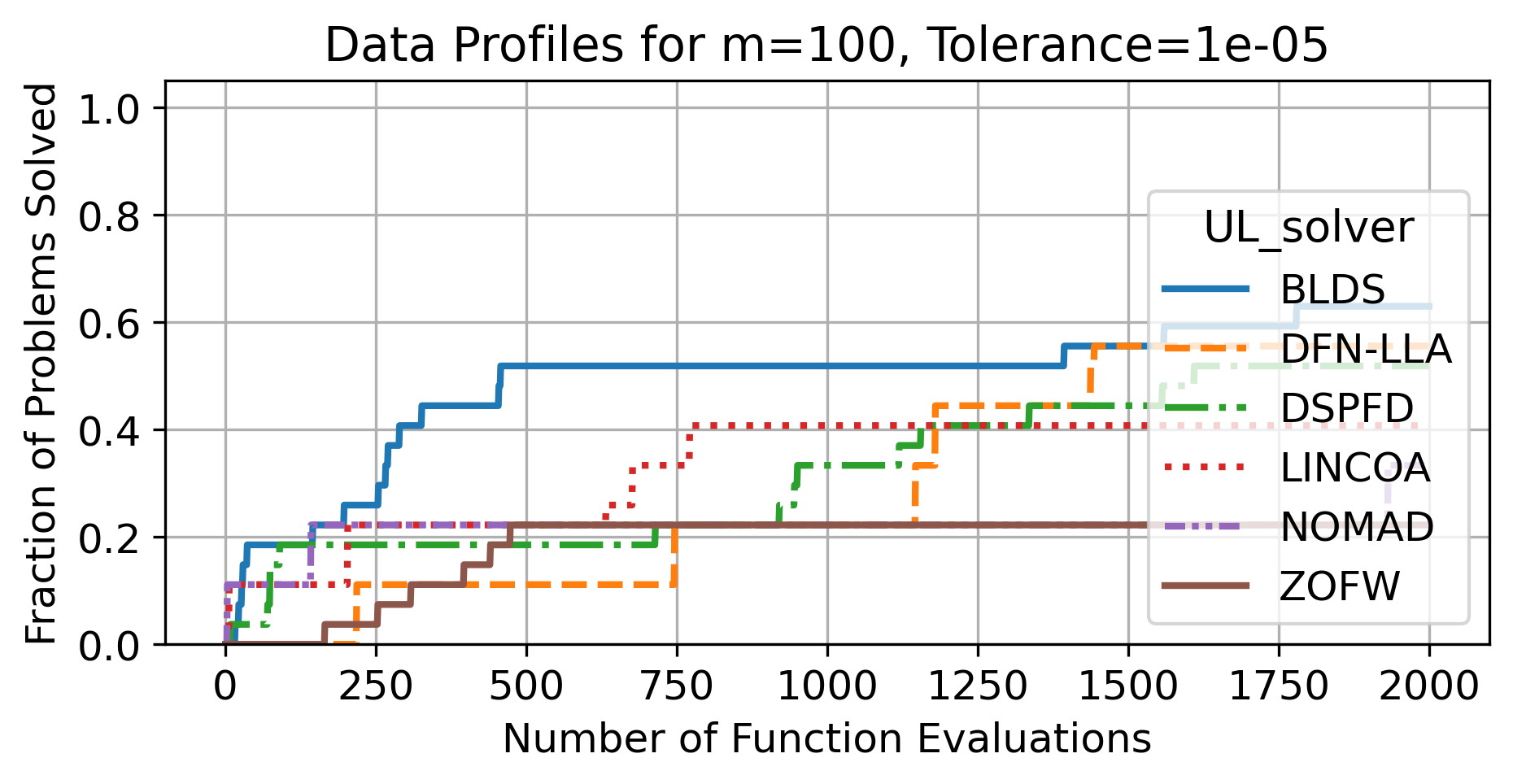}
        \caption{$m=100$}
    \end{subfigure}
    \hfill
    \begin{subfigure}[t]{0.32\linewidth}
        \centering
        \includegraphics[width=\linewidth]{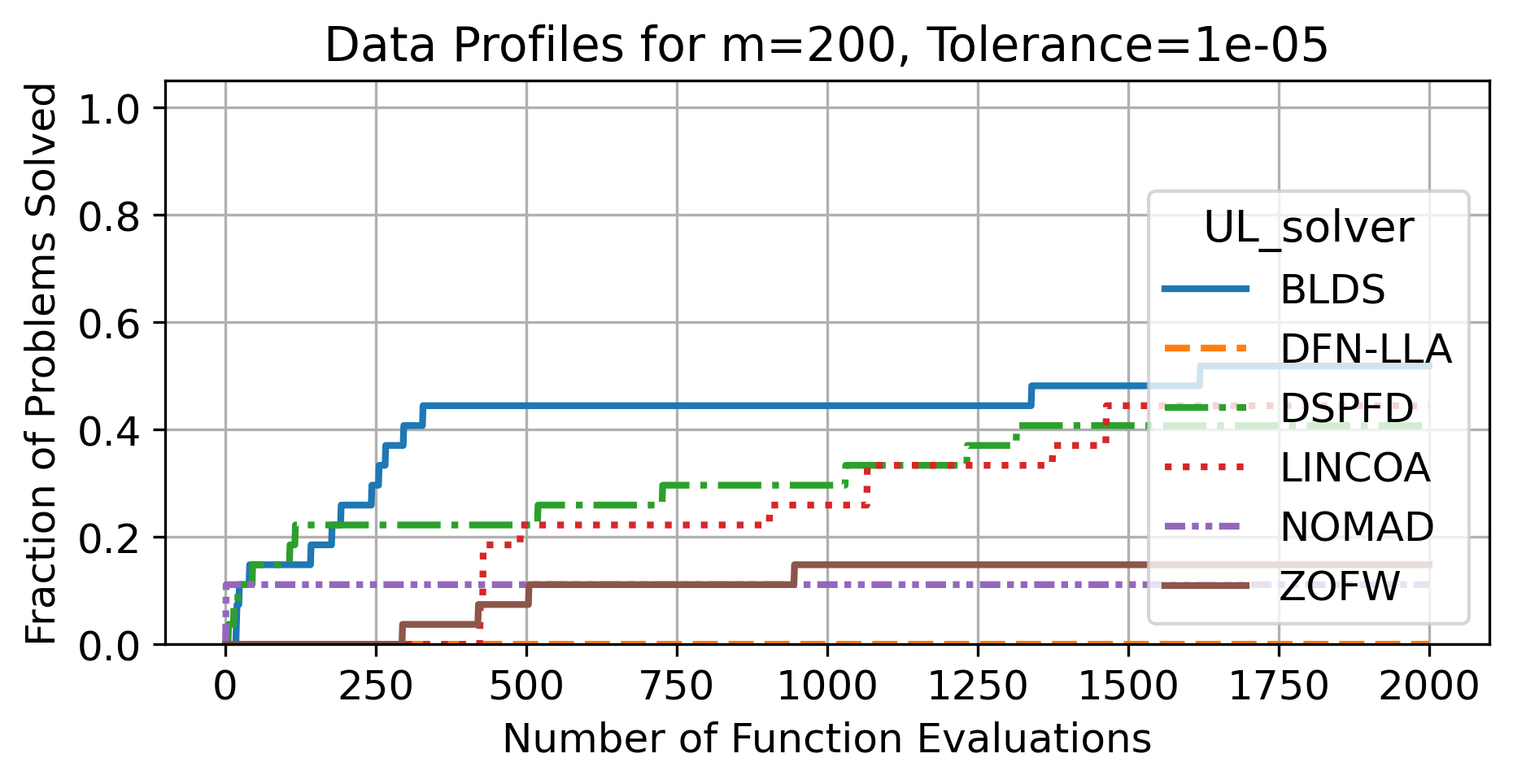}
        \caption{$m=200$}
    \end{subfigure}
    \caption{Data profiles at $\epsilon_{\mathrm{tol}}=10^{-5}$ for small vertex counts.}
    \label{fig:profiles_small_vertices_tol1e5}
\end{figure}

For small instances ($m\in\{50,100,200\}$), Figures~\ref{fig:profiles_small_vertices_tol01}--\ref{fig:profiles_small_vertices_tol1e5} show the behavior at the three tolerances.

From the plots one can clearly observe that while LINCOA and DSPFD show better performances in the case of smaller instances ($m=50$), as the problem size increases ($m\in\{100,200\}$) BLDS becomes the most efficient method in all accuracy settings. It is worth mentioning that, despite its lower efficiency, DFN-LLA shows good robustness for the case $m=100$, while the proposed BLDS is the most robust method in the case $m=200$ in all accuracy settings.

\begin{figure}[htbp]
    \centering
    \begin{subfigure}[t]{0.48\linewidth}
        \centering
        \includegraphics[width=\linewidth]{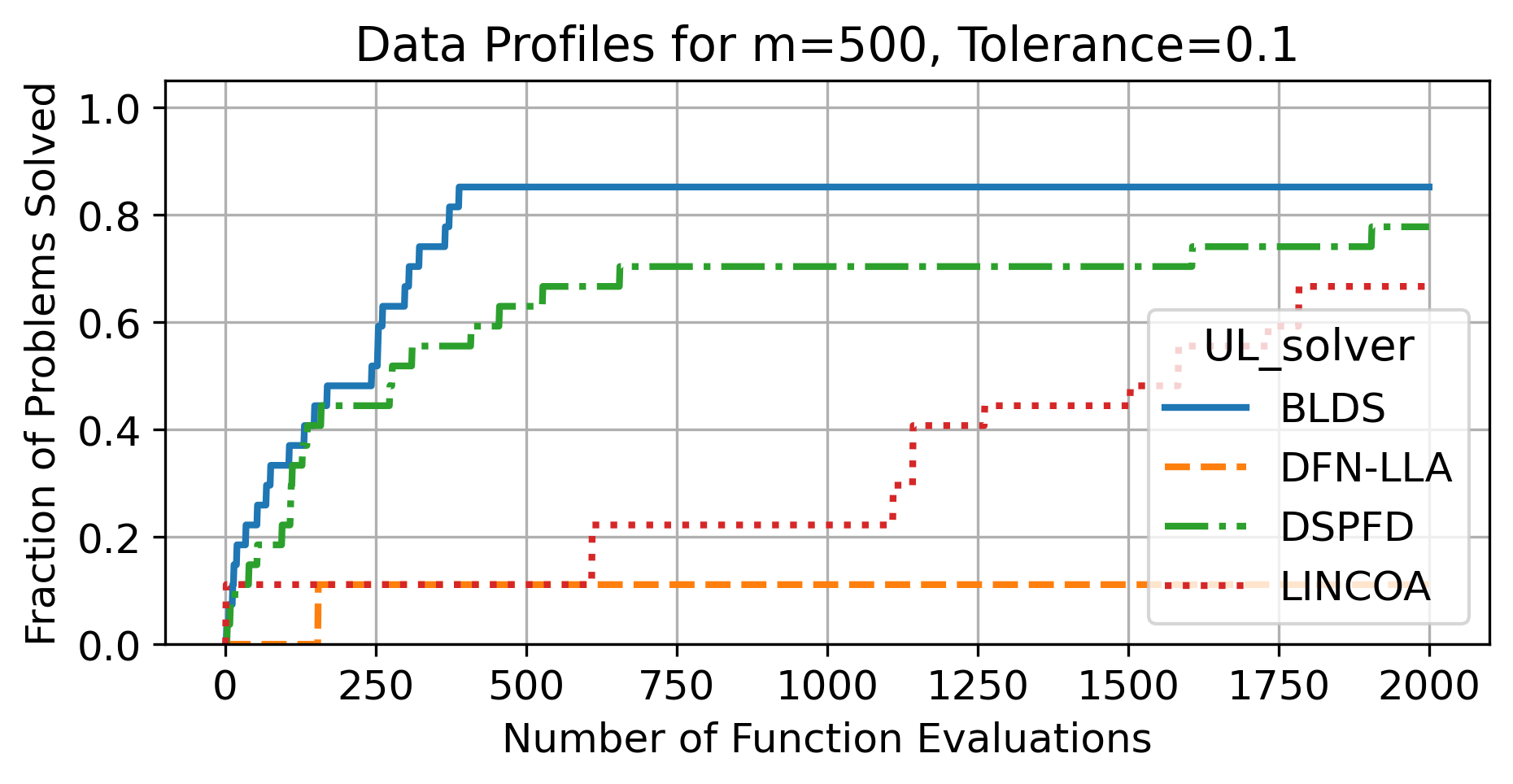}
        \caption{$m=500$}
    \end{subfigure}
    \hfill
    \begin{subfigure}[t]{0.48\linewidth}
        \centering
        \includegraphics[width=\linewidth]{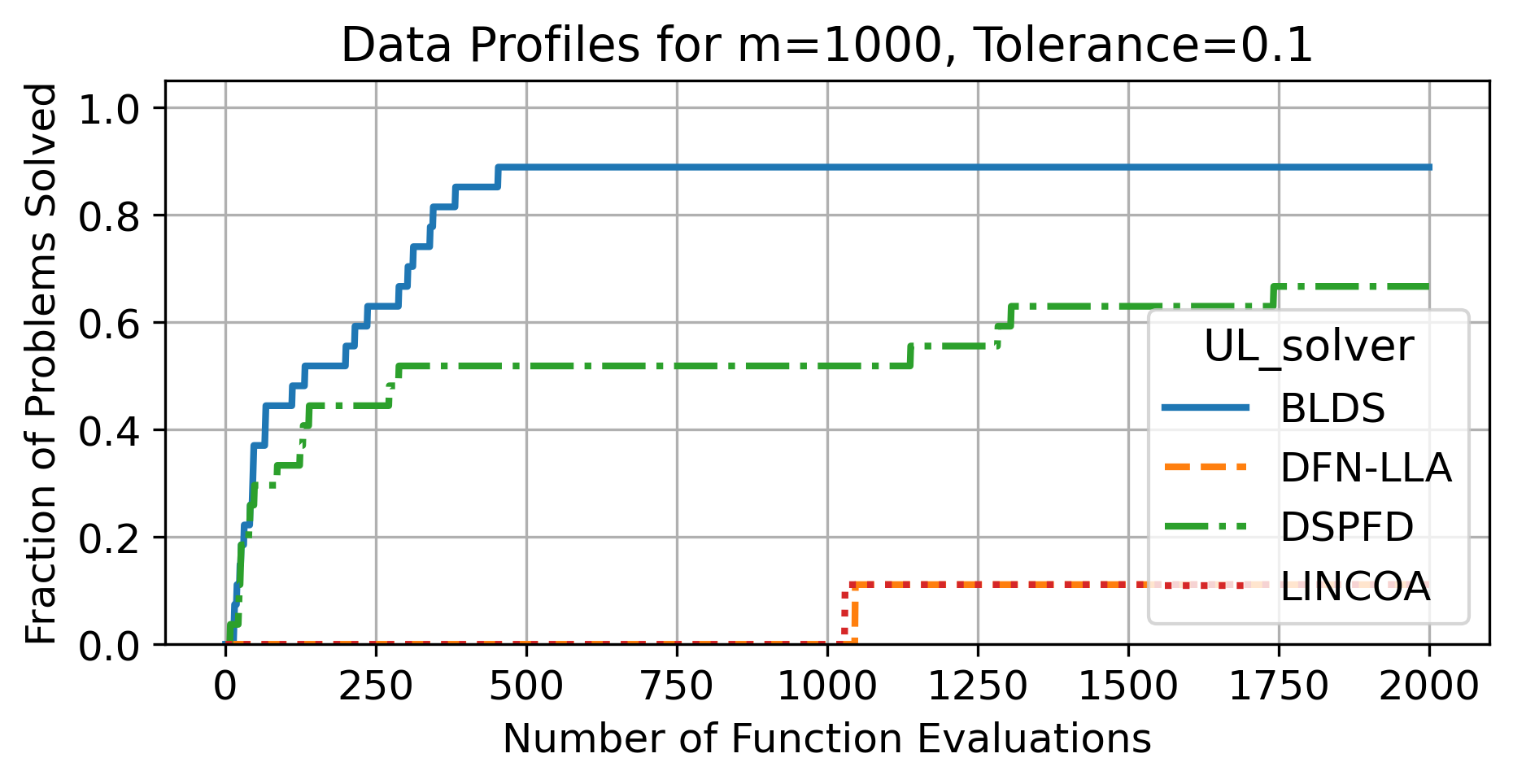}
        \caption{$m=1000$}
    \end{subfigure}
    \caption{Data profiles at $\epsilon_{\mathrm{tol}}=10^{-1}$ for large vertex counts.}
    \label{fig:profiles_large_vertices_tol01}
\end{figure}

\begin{figure}[htbp]
    \centering
    \begin{subfigure}[t]{0.48\linewidth}
        \centering
        \includegraphics[width=\linewidth]{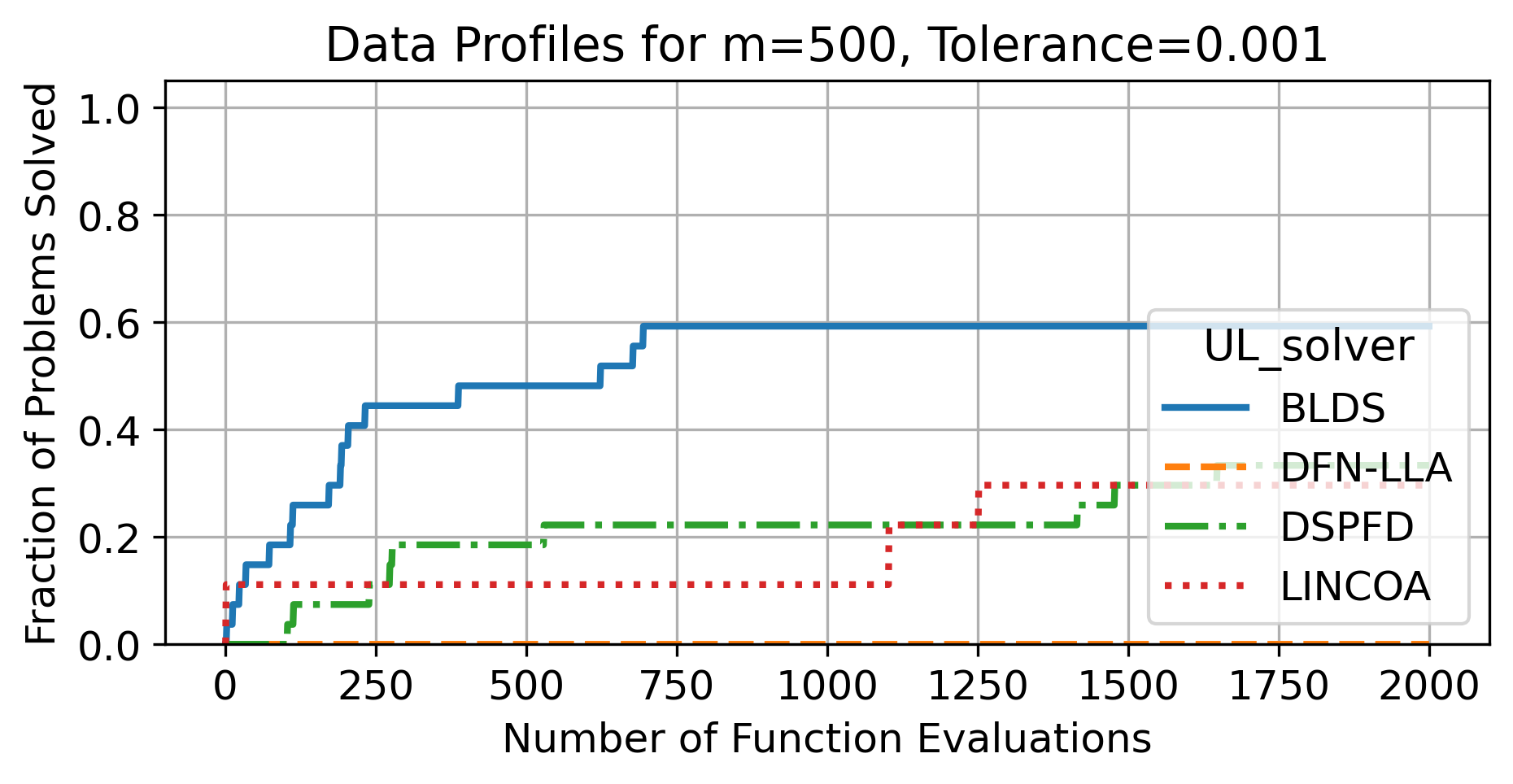}
        \caption{$m=500$}
    \end{subfigure}
    \hfill
    \begin{subfigure}[t]{0.48\linewidth}
        \centering
        \includegraphics[width=\linewidth]{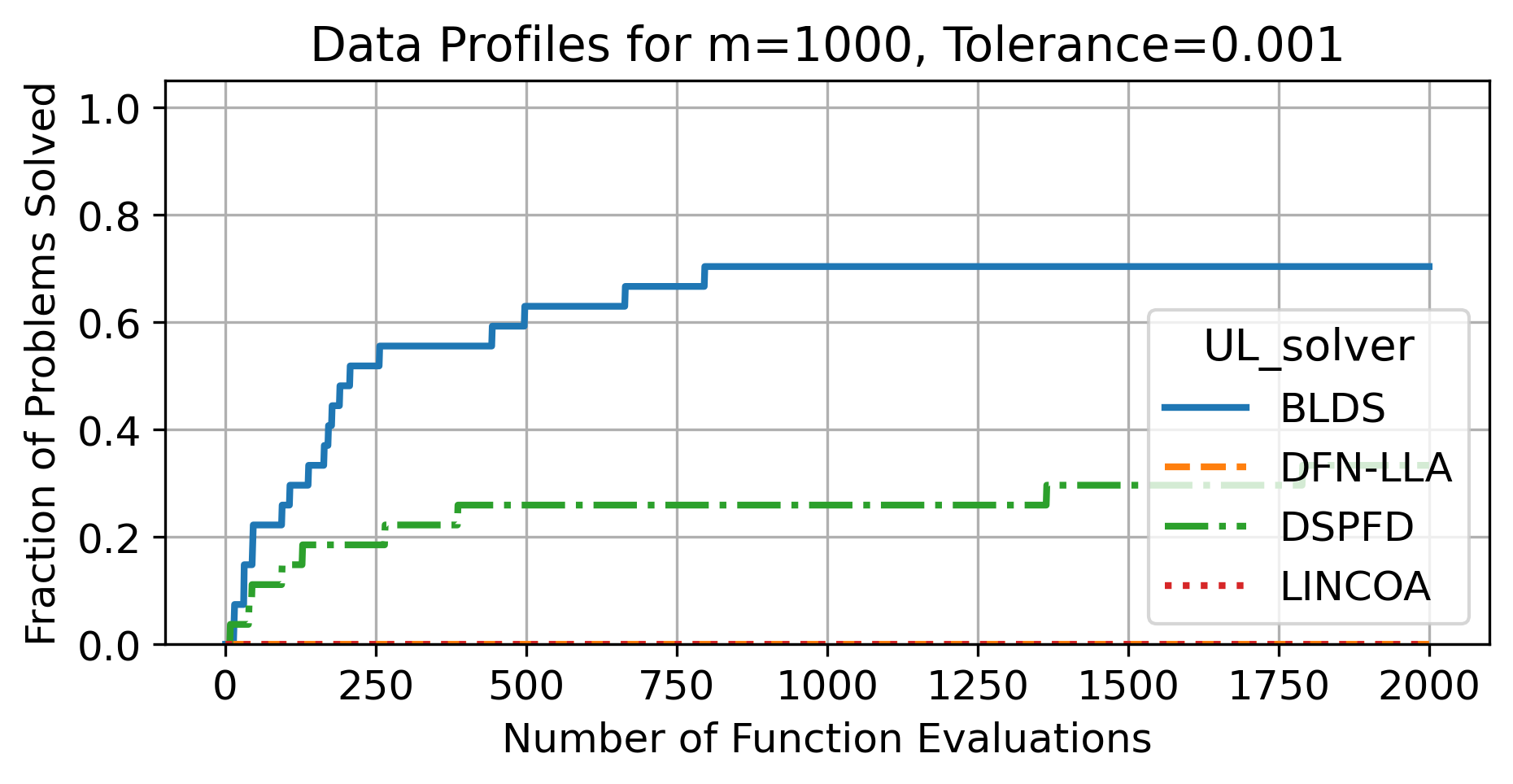}
        \caption{$m=1000$}
    \end{subfigure}
    \caption{Data profiles at $\epsilon_{\mathrm{tol}}=10^{-3}$ for large vertex counts.}
    \label{fig:profiles_large_vertices_tol001}
\end{figure}

\begin{figure}[htbp]
    \centering
    \begin{subfigure}[t]{0.48\linewidth}
        \centering
        \includegraphics[width=\linewidth]{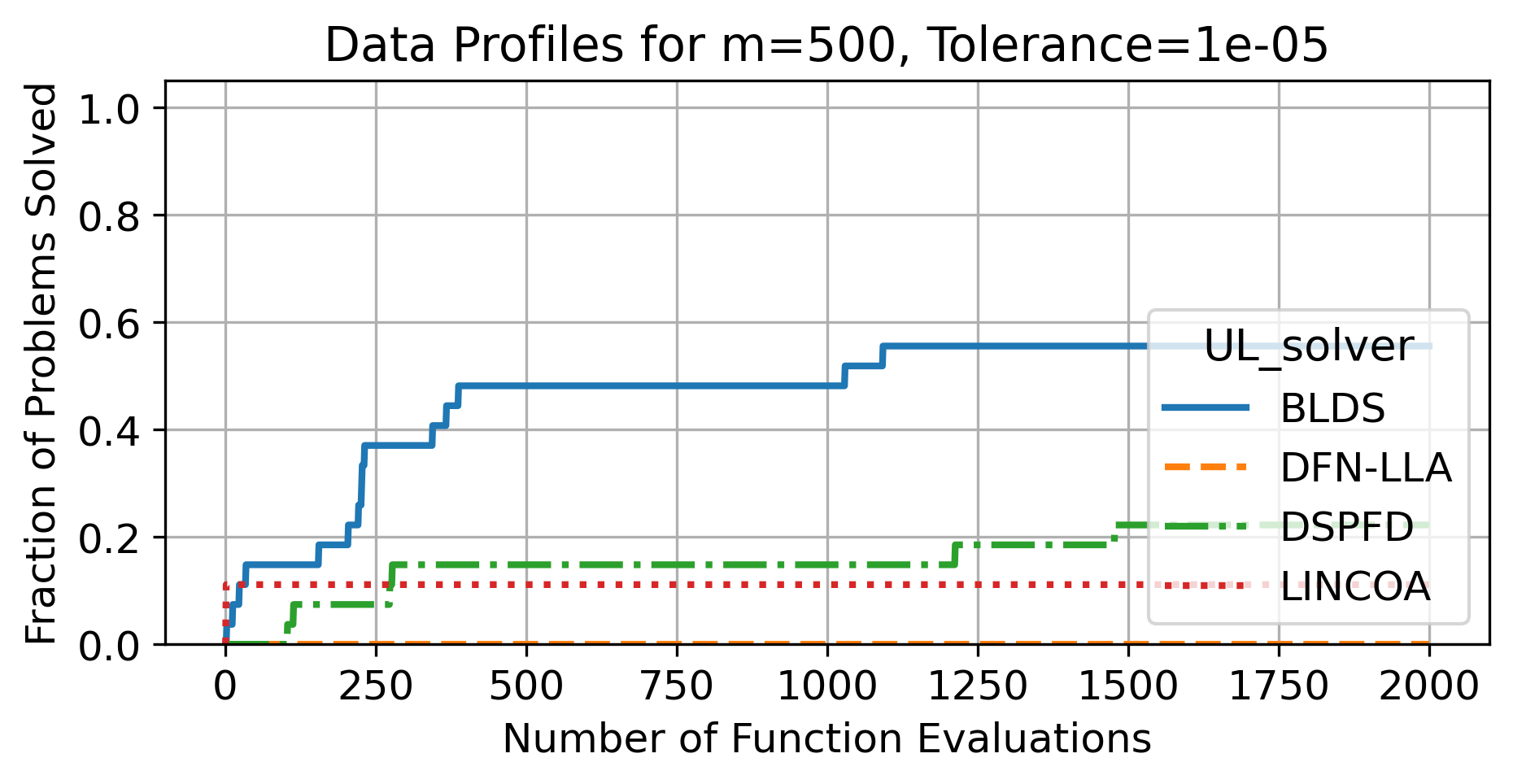}
        \caption{$m=500$}
    \end{subfigure}
    \hfill
    \begin{subfigure}[t]{0.48\linewidth}
        \centering
        \includegraphics[width=\linewidth]{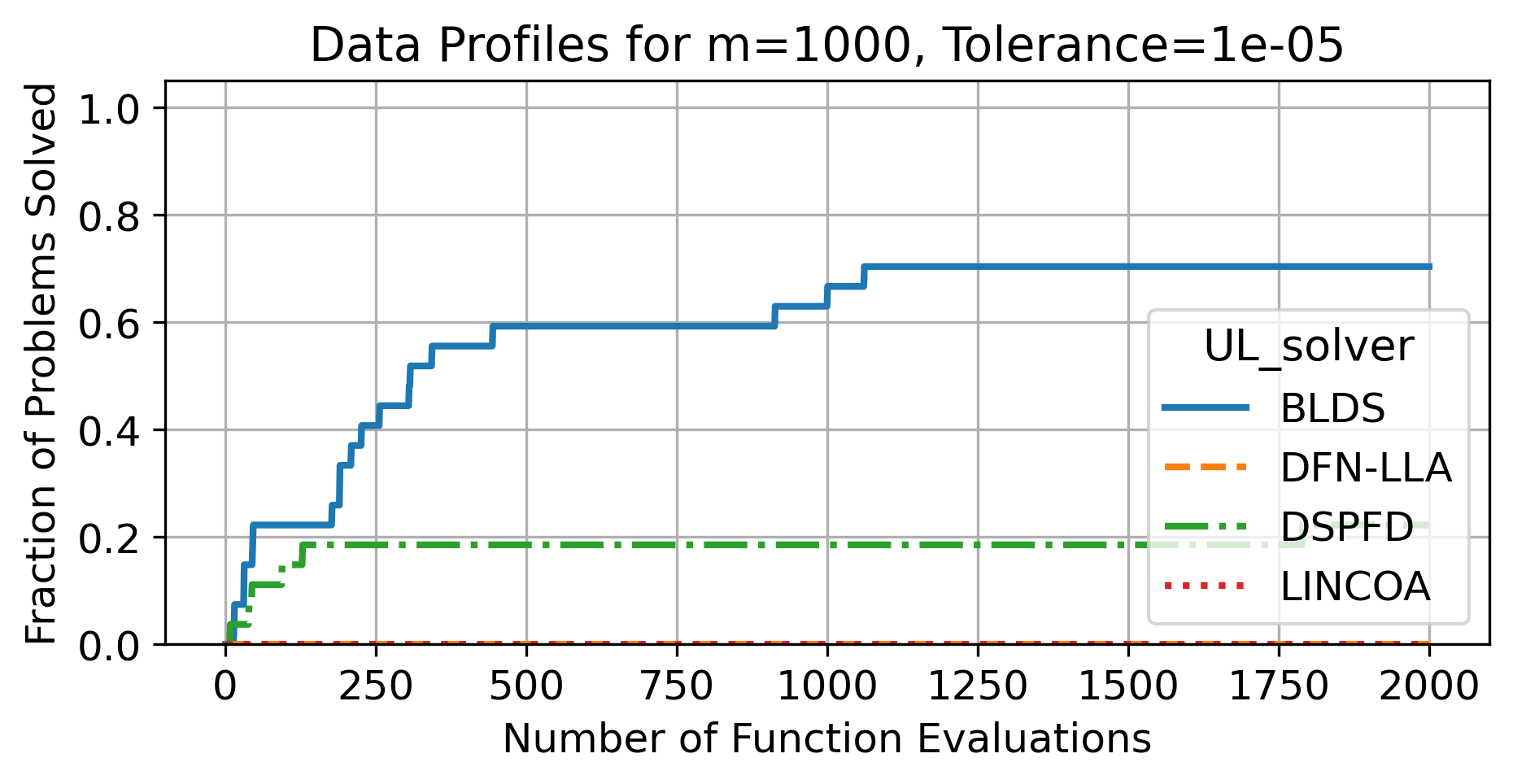}
        \caption{$m=1000$}
    \end{subfigure}
    \caption{Data profiles at $\epsilon_{\mathrm{tol}}=10^{-5}$ for large vertex counts.}
    \label{fig:profiles_large_vertices_tol1e5}
\end{figure}

In Figures~\ref{fig:profiles_large_vertices_tol01}--\ref{fig:profiles_large_vertices_tol1e5} we report the results obtained for the largest instances considered, namely $m\in\{500,1000\}$. For these plots we decided to remove NOMAD and ZOFW as they struggled with these problem sizes.

The results highlight the scalability advantage of BLDS, which for all problem sizes and accuracy settings is the most efficient and robust strategy among the four considered ones.

\section{Conclusion}
\label{sec:conclusion}

We have proposed a derivative-free direct-search algorithm for bilevel optimization problems where the upper-level feasible set is a polytope. By reformulating the problem on the unit simplex via barycentric coordinates, the method exploits Frank-Wolfe-type directions to maintain feasibility while requiring only inexact zeroth-order evaluations from a lower-level oracle.

The main contribution of our work is the theoretical analysis, which establishes convergence guarantees in two settings. In the smooth case, we proved finite termination for a positive minimum stepsize and bounded the Frank-Wolfe gap at the final iterate, showing that a suitable choice of the minimum stepsize yields an $O(\sqrt{\varepsilon})$ gap. In the nonsmooth case, we introduced the notion of a Nearby-Clarke Frank-Wolfe stationary point (NCFWSP) and established approximate stationarity guarantees for a dense-direction variant under Lipschitz continuity of the reduced objective, accounting for both the minimum stepsize and the oracle error. These results extend the analysis of inexact direct-search methods to the constrained bilevel setting through Frank-Wolfe gap-based stationarity measures. We also performed numerical experiments to assess the performance of the proposed strategy against established derivative-free solvers applied to a suitable reformulation of the bilevel problem. The results illustrate the competitiveness of the proposed strategy on the tested instances, particularly for larger numbers of vertices and more stringent target accuracies.

In future work, we aim to analyze extensions of the proposed strategy to polyhedral feasible sets specified by linear constraints, including cases where the full list of vertices is not known a priori or is expensive to compute. We also plan to investigate computationally feasible direction-sampling strategies for the dense-direction variant used in the nonsmooth analysis.

\bibliographystyle{unsrt}
\bibliography{referencesnew}

\end{document}